\documentclass[11pt,reqno]{amsart}

\newtheorem{remark}{Remark}
\numberwithin{equation}{section}
\usepackage{amsmath}
\usepackage{amssymb}
\usepackage{amsfonts}
\usepackage{graphicx}
\usepackage{subfigure}
\usepackage{stmaryrd}
\usepackage{caption}
\usepackage{import}
\usepackage{xifthen}
\usepackage{pdfpages}
\usepackage{transparent}

\usepackage{comment}
\usepackage{tikz}
\usepackage{url}
\usepackage[linktocpage]{hyperref}
\hypersetup{
    colorlinks=true,
    citecolor=blue,
    linkcolor=blue,
    filecolor=blue,
    urlcolor=blue,}
\calclayout 

\begin{document}

\title[Thermal quasi-geostrophic model on the sphere]{Thermal quasi-geostrophic model on the sphere: derivation and structure-preserving simulation}

\author{Michael Roop}
\address{Michael Roop: Department of Mathematical Sciences, Chalmers University of Technology and University of Gothenburg, 412 96 Gothenburg, Sweden}
\email{michael.roop@chalmers.se}

\author{Sagy Ephrati}
\address{Sagy Ephrati: Department of Mathematical Sciences, Chalmers University of Technology and University of Gothenburg, 412 96 Gothenburg, Sweden}
\email{sagy@chalmers.se}



\begin{abstract}
We derive the global model of thermal quasi-geostrophy on the sphere via asymptotic expansion of the thermal rotating shallow water equations. The model does not rely on the asymptotic expansion of the Coriolis force and extends the quasi-geostrophic model on the sphere by including an additional transported buoyancy field acting as a source term for the potential vorticity. We give its Hamiltonian description in terms of semidirect product Lie--Poisson brackets. The Hamiltonian formulation reveals the existence of an infinite number of conservation laws, Casimirs, parameterized by two arbitrary smooth functions. A structure-preserving discretization is provided based on Zeitlin's self-consistent matrix approximation for hydrodynamics. A Casimir-preserving time integrator is employed to numerically fully preserve the resulting finite-dimensional Lie--Poisson structure. Simulations reveal the formation of vorticity and buoyancy fronts, and large-scale structures in the buoyancy dynamics induced by the buoyancy-bathymetry interaction.
\end{abstract}

\maketitle

\section{Introduction}
The thermal rotating shallow water model (TRSW) is known to contain the basic mechanisms of the ocean and atmospheric dynamics on a planetary scale \cite{Zeitlin2018}, such as horizontal circulation of a fluid caused by a misalignment of horizontal gradients of buoyancy and bathymetry with simultaneous transport of buoyancy, and serves as a common model in geophysical fluid dynamics (GFD). The TRSW model is obtained via subsequent approximations of the 3D Euler equations for incompressible stratified inhomogeneous fluids. For small buoyancy stratification, it simplifies to the Euler-Boussinesq equations, and, after vertical averaging, to the thermal rotating Green--Naghdi equations. Finally, neglecting the non-hydrostatic pressure effects, one obtains the TRSW equations. For a detailed derivation, as well as for the stochastic versions of the mentioned models, we refer to \cite{HolmLuesink2021,HolmLuesinkPan2021}. The TRSW model incorporates thermal effects through the horizontally varying buoyancy field transported by the flow and describes the motion of a two-dimensional upper layer of the fluid on top of an inert lower layer with varying bottom topography (bathymetry). The TRSW equations were first outlined as early as in 1960s in the work \cite{BrienReid1967} and were further developed in \cite{Ripa1993,Ripa1995,Ripa1999}. 

The TRSW model contains a number of dimensionless parameters, such as the Rossby number, the Froude number, and the buoyancy stratification parameter. These parameters are small in the geophysical regime, which leads to further simplification of the TRSW model via asymptotic expansions. This yields the \textit{thermal quasi-geostrophic (TQG) equations}, which have been derived and analyzed on the $\beta$-plane \cite{HolmLuesinkPan2021,WarnefordDellar2013, beron2021nonlinear} and for which local in time unique strong solutions were proven to exist \cite{CrisanHolmLuesinkMensahPan2023}. The key mechanism behind the derivation of such planar versions of the TQG model is the expansion of the Coriolis parameter in a regular series with respect to the small Rossby number, which corresponds to a planar approximation of the Earth's surface in the neighborhood of a fixed latitude. The model therefore describes local dynamics away from equator, but does not constitute a global model on the sphere. In the present paper, we extend the TQG model to the entire sphere, by keeping the full variation of the Coriolis parameter in the thermal quasi-geostropic balance equation, and present a structure-preserving numerical integration method for the resulting system of equations.

Quasi-geostrophic (QG) models are the limiting cases of the TQG equations for the constant buoyancy. The former are typically derived in the planar approximation using the near-balance of the Coriolis force and the pressure gradient, as well as either the $f$-plane or the $\beta$-plane approximation to the Coriolis parameter at a certain latitude. For the $f$-plane approximation, the Coriolis parameter is merely replaced with the constant value $f_{0}$, whereas for the $\beta$-plane approximation, one takes the linear part of the expansion, $f=f_{0}+\beta y$, where $y$ is the vertical coordinate in the plane tangent to the Earth's surface. This approach is known as "theoretician's geostrophy" \cite{Blackburn1985} and faces a significant difficulty when attempting to extend it to the entire sphere. Namely, there is no other $f_{0}$ than zero that would represent the whole sphere. The alternative "synoptician's geostrophy" by Blackburn \cite{Blackburn1985} suggests to keep the full variation of $f$ in the geostrophic balance. However, this leads to the meridional velocity vanishing at the equator, which is also unrealistic. Perhaps the first model that overcomes the mentioned difficulties was developed by Lorenz \cite{Lorenz1960}, who, after subsequent simplifications of the velocity divergence equation, obtained the \textit{linear balance equation}. This was further simplified by Daley \cite{Daley1983} who found "the simplest form of the geostrophic relationship" by taking the trivial solution to Lorenz's balance equation.
In \cite{Verkley2009}, another balance relation was derived for the one-layer shallow water model based on Daley's balance equation. This relation contains the Cressman stretching term \cite{Cressman1958} proportional to $f^{2}$. It was further re-derived in the recent work \cite{LFEG2024} via the perturbation series in vorticity and velocity divergence for the rotating shallow water (RSW) model. The perturbation series is a promising approach in application to the TRSW model, as it allows to derive quasi-geostrophic models and specify valid parameter regimes. In the present work, we use the perturbation series to obtain the global TQG model on the sphere and show how the buoyancy field contributes to the geostrophic balance relation.

Geophysical fluid-dynamical models often possess a Hamiltonian formulation. This observation originates from the seminal work of Arnold \cite{Arn}, and the QG models are no exception. The Hamiltonian framework provides a systematic way to establish the infinite number of conservation laws in 2D hydrodynamics called \textit{Casimirs}. When it comes to computer simulations of the dynamics of the mentioned models, it is crucial to preserve the Casimir invariants to guarantee stable numerical solutions and accurate prediction of long-time statistics. A Casimir-preserving discretization was developed by \cite{ModViv,ModViv1} for the Euler equations on the sphere and was used to study long-time solution behavior. It further provided numerical evidence for the double cascade in two-dimensional turbulence \cite{CiViMo2023,ModViv2} conjectured by Kraichnan \cite{kraichnan1967inertial}. The discretization builds on Zeitlin's self-consistent truncation \cite{Zeit}, which allows approximating the infinite-dimensional Lie--Poisson structure by its finite-dimensional matrix counterpart. The global QG model possess a similar Lie--Poisson structure \cite{LFEG2024}, which has enabled structure-preserving numerical methods for single-layer \cite{franken2024singlelayer} and multi-layer \cite{franken2025casimir} QG models. Numerical results revealed the formation of stable zonal jets in long-time simulations. The mentioned models are essentially one-field transport equations, meaning that the prognostic field (the vorticity in the Euler equations, or the potential vorticity in the QG equations) is advected by the corresponding stream function related to the vorticity via the Laplace or Helmholtz operator. Their matrix approximations are isospectral flows, which makes it possible to utilize the isospectral time integrator developed in \cite{ModViv}.

A Hamiltonian formulation can also be found for the global TQG model as presented in this paper. Global TQG dynamics are described by a multi-field model with a single transported quantity, rather than a single-field transport model. Namely, the presence of varying buoyancy breaks the symmetry of the QG equations and introduces a source term in the vorticity advection equation. This leads to the loss of enstrophy as a conserved quantity, which is a Casimir in the QG model. Instead, the TQG equations are formulated in terms of the \textit{semidirect product Lie--Poisson bracket}. This bracket was originally developed in the context of magnetohydrodynamics (MHD) \cite{HMR,MoGr1980,HolmKuper1983}, where the magnetic field is transported by the fluid, while the Lorentz force is added to the vorticity equation and thereby breaks the symmetry of the original Euler equations. Thus, from a mathematical perspective, the magnetic field plays the same role in MHD as the buoyancy does in TQG. The appearance of the semidirect product bracket in TQG is somewhat expected, as it appears whenever the symmetry is broken \cite{KhMisMod}. The Hamiltonian formulation of the TQG model allows identifying Casimir invariants resembling those of MHD. Similarly, the Zeitlin matrix truncation can be found for TQG and the resulting system can be integrated in time using the \textit{magnetic midpoint integrator} developed in \cite{ModRoop} for plasma-physical models. The integrator exactly preserves the Casimir invariants and nearly preserves the energy.

The paper is organized as follows. In section~\ref{sect2}, we derive the global TQG model on the sphere via asymptotic expansion of the the TRSW equations with respect to a small parameter that unifies the Froude number, the Rossby number, and the stratification parameter. In section~\ref{sect3}, we give a Hamiltonian formulation of the TQG model in terms of a non-canonical Hamiltonian structure, along with the Casimir invariants. In section~\ref{sect4}, we present the structure-preserving discretization for TQG on the sphere, and demonstrate simulation results in section~\ref{sect5}. We conclude the paper in section ~\ref{sect6}.

\vspace{0.5cm}
\noindent
{\bf Acknowledgements.} This work was supported by the Knut and Alice Wallenberg Foundation, grant number WAF2019.0201, and by the Swedish Research Council, grant number 2022-03453. The work of M.R. is supported by the grants from the Royal Swedish Academy of Sciences (MA2024-0034, MG2024-0050). The authors would like to thank Erwin Luesink, Arnout Franken, and Darryl Holm for inspiring discussions.

\section{Model derivation}
\label{sect2}
\begin{figure}[ht!]
\centering
\includegraphics[scale=0.7]{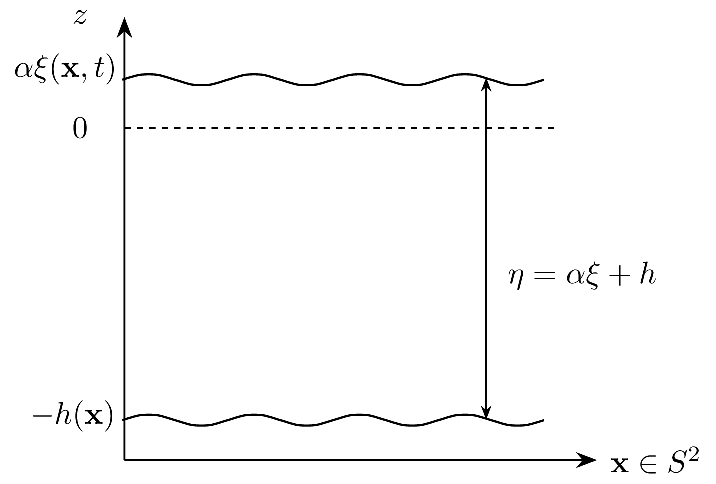}
\caption{Vertical structure of the flow domain. The free surface is given by the function $\alpha\xi(\mathbf{x},t)$, which is small compared to the total depth $\eta=\alpha\xi+h$, with $h(\mathbf{x})$ being the bathymetry.}
\label{layer}
\end{figure}

In this section, we derive the global TQG model on the sphere starting from the TRSW equations by following the approach presented in \cite{LFEG2024}. In what follows, we denote the sphere by $S^{2}\subset\mathbb{R}^{3}$ and let $\mathbf{x}\in S^{2}$.
The non-dimensional TRSW equations read \cite{HolmLuesinkPan2021}
\begin{equation}
\label{TRSW}
\left\{
\begin{aligned}
\displaystyle&\frac{\partial\mathbf{u}}{\partial t}+(\mathbf{u}\cdot\nabla)\mathbf{u}+\frac{1}{\mathrm{Ro}}f\mathbf{z}\times\mathbf{u}=-\frac{\alpha}{\mathrm{Fr}^{2}}\nabla((1+sb)\xi)+\frac{s}{2\mathrm{Fr}^{2}}(\alpha\xi-h)\nabla b, \\
\\
\displaystyle&\frac{\partial\eta}{\partial t}+\nabla\cdot(\eta\mathbf{u})=0,\quad\frac{\partial b}{\partial t}+(\mathbf{u}\cdot\nabla)b=0,
\end{aligned}
\right.
\end{equation}
where $\mathbf{u}(\mathbf{x},t)$ is the velocity field of the fluid (note that at this point there is no assumption of vanishing divergence), $\mathbf{z}$ is the outward unit normal vector on the sphere, $f=2\cos(\theta)$ is the dimensionless Coriolis parameter, with $\theta$ being the latitude ($\theta=\pi/2$ at the equator); $\alpha\xi(\mathbf{x},t)$ is the free surface elevation, so that $\eta(\mathbf{x},t)=\alpha\xi(\mathbf{x},t)+h(\mathbf{x})$ is the total depth, and $h(\mathbf{x})=1+\varepsilon h_{1}(\mathbf{x})$ is the bathymetry function. The constant $\alpha$ is the typical wave amplitude, which is assumed to be small compared to the total depth $\eta(\mathbf{x},t)$. We illustrate this in Fig.~\ref{layer}.

The field $b(\mathbf{x},t)=(\rho(\mathbf{x},t)-\rho_{0})/\rho_{0}$ is the dimensionless buoyancy defined as a normalized fluid density variation. The system also contains three dimensionless parameters. The Rossby number $\mathrm{Ro}=U/(\Omega L)$ is the ratio between the typical velocity $U$ and the rotation velocity $\Omega L$ ($\Omega$ is the Earth's rotation frequency), with $L$ being the typical length scale. The stratification parameter $s$ governs the importance of buoyancy. The Froude number $\mathrm{Fr}=U/\sqrt{gH}$ is the ratio between the typical velocity $U$ and the speed of the fastest gravity wave $\sqrt{gH}$, with $g$ being the gravitational acceleration and $H$ being the average fluid depth. In the geophysical approximation, the parameters $\mathrm{Ro}$, $\mathrm{Fr}$, $\alpha$, $s$ are assumed to be small of the same order, and for convenience we shall denote them by the same symbol $\varepsilon$:
\begin{equation*}
\mathcal{O}(\alpha)=\mathcal{O}(s)=\mathcal{O}(\mathrm{Ro})=\mathcal{O}(\mathrm{Fr})=\mathcal{O}(\varepsilon).
\end{equation*}

\subsection{Simplest TQG balance equation}
We first derive the expression for thermal geostrophic balance. This amounts to expanding the perturbation series of several variables in \eqref{TRSW} and finding the velocity field that ensures a balance at leading order.

We start with the transport equation for $\eta(\mathbf{x},t)$. Combining the above expressions for $h(\mathbf{x},t)$ and $\eta(\mathbf{x},t)$, we get that $\eta(\mathbf{x},t)=1+\alpha\xi(\mathbf{x},t)+\varepsilon h_{1}(\mathbf{x})$. Plugging this expression to the transport equation for $\eta(\mathbf{x},t)$ in \eqref{TRSW}, we obtain
\begin{equation}
\label{etaeqexpansion}
\alpha\frac{\partial\xi}{\partial t}+\nabla(\alpha\xi+\varepsilon h_{1})\cdot\mathbf{u}+\nabla\cdot\mathbf{u}+(\alpha\xi+\varepsilon h_{1})\nabla\cdot\mathbf{u}=0.
\end{equation}
Since $\alpha$ and $\varepsilon$ are by assumption small parameters of the same order, we get from \eqref{etaeqexpansion} that the leading term in the velocity field expansion is divergence-free. 
Equivalently, the divergence of the velocity field $\mathbf{u}$ is of order $\varepsilon$. Small horizontal divergence has been a common assumption in quasi-geostrophic theories and has recently been validated on observational data \cite{YanoMuletBonazzola2009}.
Combining this fact with the Helmholtz decomposition of vector fields on the sphere, we obtain
\begin{equation}
\label{veolcitydec}
\mathbf{u}=\mathbf{z}\times\nabla\psi+\varepsilon\nabla\chi,
\end{equation}
where $\psi(\mathbf{x},t)$ is the \textit{stream function} generating the divergence-free part of the velocity field, and $\chi(\mathbf{x},t)$ is the potential for the gradient part of the vector field $\mathbf{u}(\mathbf{x},t)$. This observation allows us to get the leading term of the buoyancy transport equation:
\begin{equation*}
\dot b=-(\mathbf{z}\times\nabla\psi)\cdot \nabla b=\left\{b,\psi\right\},
\end{equation*}
where $\left\{\cdot,\cdot\right\}$ is the Poisson bracket on $C^{\infty}(S^{2})$.

We view the two-dimensional velocity field $\mathbf{u}(\mathbf{x},t)$ as a three-dimensional vector field with a trivial radial component, which allows for the use of the standard $\mathrm{curl}$ operator in three dimensions. Then, using the identity
\begin{equation*}
\mathrm{curl}(\mathbf{u})\times\mathbf{u}=(\mathbf{u}\cdot\nabla)\mathbf{u}-\frac{1}{2}\nabla|\mathbf{u}|^{2},
\end{equation*}
we can rewrite the velocity equation in \eqref{TRSW} as follows:
\begin{equation}
\label{velTRSW}
\frac{\partial\mathbf{u}}{\partial t}+\left(\omega+\frac{f}{\mathrm{Ro}}\right)\mathbf{z}\times\mathbf{u}=-\frac{\alpha}{\mathrm{Fr}^{2}}\nabla\left((1+sb)\xi\right)-\frac{1}{2}\nabla|\mathbf{u}|^{2}+\frac{s}{2\mathrm{Fr}^{2}}(\alpha\xi-h)\nabla b,
\end{equation}
where $\omega(\mathbf{x},t)$ is the vorticity function defined as $\omega\,\mathbf{z}=\mathrm{curl}(\mathbf{u})$.

Let us now insert the velocity decomposition \eqref{veolcitydec} in \eqref{velTRSW} and apply the divergence operator to the both sides. We get
\begin{equation*}
\varepsilon\Delta\dot\chi+\operatorname{div}\left(-\left(\omega+\frac{f}{\mathrm{Ro}}\right)\nabla\psi+\frac{\alpha}{\mathrm{Fr}^{2}}\nabla((1+sb)\xi)+\frac{1}{2}\nabla|\nabla\psi|^{2}-\frac{s}{2\mathrm{Fr}^{2}}(\alpha\xi-h)\nabla b\right)+\mathcal{O}(\varepsilon)=0,
\end{equation*}
where we used that $\mathbf{z}\times\mathbf{u}=-\nabla\psi+\mathcal{O}(\varepsilon)$. Recall that $\alpha, s, \mathrm{Ro}$, and $\mathrm{Fr}$ are all of order $\varepsilon$. Then, by collecting the terms of order $\mathcal{O}(1/\varepsilon)$, we find 
\begin{equation}
\label{TQGbal}
\nabla\cdot\left(-\frac{f}{\mathrm{Ro}}\nabla\psi+\frac{\alpha}{\mathrm{Fr}^{2}}\nabla\xi+\frac{s}{2\mathrm{Fr}^{2}}h\nabla b\right)=0.
\end{equation}

We now have to find a solution to equation \eqref{TQGbal}. A similar equation, the \textit{linear balance equation}, was derived in \cite{Lorenz1960}. A conventional simplifying assumption in GFD is that the the Coriolis parameter $f$ has a small variation, i.e. $\nabla f\sim \mathcal{O}(\varepsilon)$, which would imply that $f\nabla\psi\approx\nabla(f\psi)$ (up to $\mathcal{O}(\varepsilon)$ terms). This assumption is evident for the $f$-plane and $\beta$-plane approximations. To the best of the authors' knowledge, no rigorous justification of this fact has been offered for the spherical case. Our subsequent derivation of the TQG balance relies on the assumption that $f\nabla\psi\approx\nabla(f\psi)$. We motivate this assumption as follows. First, we recall that equations \eqref{TRSW} are dimensionless and therefore differential operators, such as the gradient, are dimensionless as well. We thus compare the two terms $f\nabla\psi$ and $\psi\nabla f$,
\begin{equation*}
f\nabla\psi=\frac{2L}{a}\cos(\theta)\frac{\partial\psi}{\partial\theta}\mathbf{e}_{\theta}+\frac{2L}{a}\cot(\theta)\frac{\partial\psi}{\partial\varphi}\mathbf{e}_{\varphi},\quad\psi\nabla f=-\frac{2L}{a}\psi\sin(\theta)\mathbf{e}_{\theta},
\end{equation*}
where $\mathbf{e}_{\theta}$ and $\mathbf{e}_{\varphi}$ are the unit vectors in polar and azimuthal directions respectively, and $a$ is the radius of the Earth. We observe that the term $\psi\nabla f$ is bounded for all $\theta\in[0,\pi]$ provided $\psi$ is bounded. Furthermore, $\psi\nabla f$ is of order $\mathcal{O}(\varepsilon)$ when $L/a \sim\mathcal{O}(\varepsilon)$, i.e., when the horizontal lengthscale is small compared to the Earth's radius. 
At the same time, despite the presence of the same multiplier $L/a$ in the term $f\nabla\psi$, the $\cot(\theta)$ function has a singularity and therefore cannot be neglected.

The simplifying assumption $f\nabla\psi\approx\nabla(f\psi)$ was used in \cite{Verkley2009}, and its validity was illustrated on the dynamics of linear Rossby waves. The global QG model on the sphere derived by means of the above assumption was investigated also in \cite{SchubTaftSilv2008}. In the work \cite{LFEG2024}, this assumption facilitates the study of global QG dynamics and reveals the formation of zonal jets.

Using the above reasoning and also that $\nabla h\sim \mathcal{O}(\varepsilon)$ we simplify equation \eqref{TQGbal} further to obtain
\begin{equation*}
\Delta\left(-f\psi+\frac{\alpha\mathrm{Bu}}{\mathrm{Ro}}\xi+\frac{s\mathrm{Bu}}{2\mathrm{Ro}}hb\right)=0,
\end{equation*}
where we introduced the Burger number as $\mathrm{Bu}=(\mathrm{Ro}/\mathrm{Fr})^{2}=O(1)$. The simplest solution is trivial and gives rise to \textit{the simplest form of TQG balance:}
\begin{equation*}
f\psi=\frac{\alpha\mathrm{Bu}}{\mathrm{Ro}}\xi+\frac{s\mathrm{Bu}}{2\mathrm{Ro}}hb.
\end{equation*}
Further, using the asymptotic expansion for the depth in dimensionless form $h(\mathbf{x})=1+\varepsilon h_{1}(\mathbf{x})$ and keeping the terms of order $O(1)$, we obtain the final expression for TQG balance:
\begin{equation}
\label{TQGbalancefinal}
f\psi=\frac{\alpha\mathrm{Bu}}{\mathrm{Ro}}\xi+\frac{s\mathrm{Bu}}{2\mathrm{Ro}}b.
\end{equation}
\begin{remark}
Equation \eqref{TQGbalancefinal} generalizes existing results in two aspects. 

Firstly, the TQG balance equation found here is the spherical extension of the balance derived in \cite{HolmLuesinkPan2021} for the planar TQG model. Indeed, in \cite{HolmLuesinkPan2021} one has 
\begin{equation*}
\psi=\xi+\frac{b}{2},
\end{equation*}
which is recovered (up to multipliers containing the dimensionless parameters) from \eqref{TQGbalancefinal} upon expanding the Coriolis parameter as $f=1+\varepsilon f_{1}$ ($\beta$-plane approximation), which is a valid assumption if one approximates the Earth's surface as a plane in the neighborhood of a certain latitude sufficiently far from the equator.

Secondly, \eqref{TQGbalancefinal} can be regarded as the thermal generalization of the spherical QG balance obtained in \cite{LFEG2024,Verkley2009}, by including buoyancy effects. Indeed, the QG balance equation derived in \cite{LFEG2024} is reconstructed from \eqref{TQGbalancefinal} by setting $b(\mathbf{x},t)=0$.
\end{remark}
\subsection{Potential vorticity equation}
We introduce the potential vorticity $q_{trsw}(\mathbf{x},t)$ for the TRSW equations,
\begin{equation}
\label{QRSW}
q_{trsw}=\frac{1}{\eta}\left(\omega+\frac{f}{\mathrm{Ro}}\right),
\end{equation}
and proceed to derive the TQG potential vorticity by expanding $q_{trsw}$ using the perturbation series.

Applying the $\mathrm{curl}$ operator to equation \eqref{veolcitydec}, and using the definition of the vorticity function $\omega\,\mathbf{z}=\mathrm{curl}(\mathbf{u})$, we get
\begin{equation*}
\omega\,\mathbf{z}=\mathrm{curl}(\mathbf{u})=\mathrm{curl}(\mathbf{z}\times\nabla\psi)+\mathcal{O}(\varepsilon)=(\Delta\psi)\mathbf{z}+\mathcal{O}(\varepsilon),
\end{equation*}
which gives the expansion for the vorticity
\begin{equation*}
\omega=\Delta\psi+\mathcal{O}(\varepsilon).
\end{equation*}
The approximation for $\eta(t,x)$ uses $h(\mathbf{x})=1+\varepsilon h_{1}(\mathbf{x})+\mathcal{O}(\varepsilon^{2})$ and reads
\begin{equation}
\label{etaappr}
\eta=\alpha\xi+h=1+\varepsilon h_{1}+\alpha\xi\Longrightarrow\frac{1}{\eta}=1-\varepsilon h_{1}-\alpha\xi+\mathcal{O}(\varepsilon^{2}).
\end{equation}
Substituting \eqref{etaappr} into \eqref{QRSW}, we obtain
\begin{equation}
\label{qrswappr}
q_{trsw}=(1-\varepsilon h_{1}-\alpha\xi)\left(\omega+\frac{f}{\mathrm{Ro}}\right)=\Delta\psi+\frac{f}{\mathrm{Ro}}-\frac{\varepsilon}{\mathrm{Ro}} h_{1}f-\frac{\alpha}{\mathrm{Ro}}\xi f+\mathcal{O}(\varepsilon)=q+\mathcal{O}(\varepsilon),
\end{equation}
where 
\begin{equation}
\label{qprefinal}
q=\Delta\psi+\frac{f}{\mathrm{Ro}}-\frac{\varepsilon}{\mathrm{Ro}}h_{1}f-\frac{\alpha}{\mathrm{Ro}}\xi f
\end{equation}
is the leading term of $q_{trsw}$.
The field $q$ is the TQG potential vorticity and will henceforth simply be referred to as the potential vorticity.

Using the TQG balance equation \eqref{TQGbalancefinal}, we express the surface elevation function $\xi$ in terms of $f,\psi,b$ and insert it into \eqref{qprefinal} to obtain the relation between the potential vorticity $q$ and the stream function $\psi$. We find
\begin{equation*}
q=\Delta\psi+\frac{f}{\mathrm{Ro}}-\frac{1}{\mathrm{Bu}}f^{2}\psi+\frac{s}{2\mathrm{Ro}}fb-\frac{\varepsilon}{\mathrm{Ro}}h_{1}f,
\end{equation*}
which can be rewritten, using Lamb's parameter $\gamma=4/\mathrm{Bu}$ and $\mu=\cos(\theta)$, to
\begin{equation}
\label{QPSieq}
q=(\Delta-\gamma\mu^{2})\psi+\frac{2\mu}{\mathrm{Ro}}-\frac{2\varepsilon}{\mathrm{Ro}}\mu h_{1}+\frac{s}{\mathrm{Ro}}\mu b.
\end{equation}
\begin{remark}
Relation \eqref{QPSieq} does not rely on the asymptotic expansion of the Coriolis parameter $f$ and only uses the simplifying assumption $f\nabla\psi\approx\nabla(f\psi)$. One can compare this relation to the one obtained in \cite{HolmLuesinkPan2021} for the $\beta$-plane approximation by assuming $f=1+\varepsilon f_{1}$,
\begin{equation}
\label{vortstreamDarryl}
q=(\Delta-1)\psi+f_{1}.
\end{equation}
An evident distinction between the relation \eqref{QPSieq} on the sphere and \eqref{vortstreamDarryl} on the $\beta$-plane is that the former contains an inhomogeneous (latitude-dependent) Helmholtz operator, whereas the latter is described by a homogeneous Helmholtz operator.

Furthermore, one can compare \eqref{QPSieq} to the similar relation between $q$ and $\psi$ derived in \cite{LFEG2024}, which is reconstructed (up to the sign convention) from \eqref{QPSieq} by setting $b(\mathbf{x},t)=0$.
\end{remark}

Finally, we derive the evolution equation of the potential vorticity $q$. First, let us apply the $\mathrm{curl}$ operator to the momentum equation in \eqref{TRSW}, which gives 
\begin{equation}
\label{vortTRSW}
\dot\omega+\mathrm{div}(\omega\mathbf{u})+\frac{1}{\mathrm{Ro}}\mathbf{z}\cdot\mathrm{curl}(f\mathbf{z}\times\mathbf{u})=\frac{s}{2\mathrm{Fr}^{2}}\mathbf{z}\cdot(\nabla(\alpha\xi-h)\times\nabla b).
\end{equation}
Further, using the definition \eqref{QRSW} of $q_{trsw}$ and the evolution equation for $\eta$ in \eqref{TRSW}, we obtain
\begin{equation*}
\dot\omega=\eta\dot{q}_{trsw}-\mathrm{div}(\eta\mathbf{u})q_{trsw}.
\end{equation*}
Substituting this in \eqref{vortTRSW} and eliminating the vorticity $\omega$ by means of \eqref{QRSW} we find
\begin{equation}
\label{qtrswevol}
\eta\dot{q}_{trsw}-\mathrm{div}(\eta\mathbf{u})q_{trsw}+\mathrm{div}(\eta q_{trsw}\mathbf{u})-\frac{1}{\mathrm{Ro}}\mathrm{div}(f\mathbf{u})+\frac{1}{\mathrm{Ro}}\mathbf{z}\cdot\mathrm{curl}(f\mathbf{z}\times\mathbf{u})=\frac{s}{2\mathrm{Fr}^{2}}\mathbf{z}\cdot(\nabla(\alpha\xi-h)\times\nabla b).
\end{equation}
It is readily checked that the two terms containing $\mathrm{Ro}$ cancel out, and we come to the evolution equation for $q_{trsw}$, see also \cite{HolmLuesinkPan2021}:
\begin{equation}
\label{qtrswevol1}
\dot q_{trsw}+\mathbf{u}\cdot\nabla q_{trsw}=\frac{s}{2\eta\mathrm{Fr}^{2}}\mathbf{z}\cdot\left(\nabla(\alpha\xi-h)\times\nabla b\right).
\end{equation}
Using the expansions
\begin{equation*}
q_{trsw}=q+\mathcal{O}(\varepsilon),\quad\mathbf{u}=\nabla^{\perp}\psi+\mathcal{O}(\varepsilon),\quad h=1+\varepsilon h_{1}+O(\varepsilon^{2}),\quad \frac{1}{\eta}=1-\varepsilon h_{1}-\alpha\xi+O(\varepsilon^{2}),
\end{equation*}
the identity $\mathbf{z}\cdot(\nabla K\times\nabla L)=\left\{K,L\right\}$ for arbitrary $K$ and $L$, and truncating equation \eqref{qtrswevol1} at $\mathcal{O}(1)$, we find the evolution equation for the potential vorticity $q$:
\begin{equation*}
\dot q=\left\{q,\psi\right\}+\frac{s}{2\mathrm{Fr}^{2}}\left\{b,\varepsilon h_{1}-\alpha\xi\right\}.
\end{equation*}
From the TQG balance equation \eqref{TQGbalancefinal} we have
\begin{equation*}
-\alpha\xi=\frac{sb}{2}-\frac{\mathrm{Ro}}{\mathrm{Bu}}f\psi,
\end{equation*}
and the potential vorticity evolution thus reads
\begin{equation*}
\dot q=\left\{q,\psi\right\}+\left\{b,\frac{s\varepsilon}{2\mathrm{Fr}^{2}}h_{1}-\frac{s}{2\mathrm{Ro}}f\psi\right\}.
\end{equation*}

Summarizing, we arrive at the closed model of thermal quasi-geostrophy on the sphere given by
\begin{equation}
\label{TQGsphere}
\left\{
\begin{aligned}
\displaystyle&\dot q=\left\{q,\psi\right\}+\left\{b,j\right\},\\
\displaystyle&\dot b=\left\{b,\psi\right\},\\
\displaystyle&q=(\Delta-\gamma\mu^{2})\psi+\frac{2\mu}{\mathrm{Ro}}-\frac{2\varepsilon}{\mathrm{Ro}}\mu h_{1}+\frac{s}{\mathrm{Ro}}\mu b,\\
\displaystyle&j=\frac{s\varepsilon}{2\mathrm{Fr}^{2}}h_{1}-\frac{s}{\mathrm{Ro}}\mu\psi,
\end{aligned}
\right.
\end{equation}
with $\mu=\cos(\theta)$ and $\gamma$ the Lamb parameter.

For simplicity, we choose $s=\mathrm{Ro}$, $\varepsilon=\mathrm{Ro}/2$, $\mathrm{Fr}=\mathrm{Ro}/2$, such that \eqref{TQGsphere} simplifies to a system containing only the Lamb parameter $\gamma$ and the Rossby number $\mathrm{Ro}$:
\begin{equation}
\label{TQGsphere1}
\left\{
\begin{aligned}
\displaystyle&\dot q=\left\{q,\psi\right\}+\left\{b,j\right\},\\
\displaystyle&\dot b=\left\{b,\psi\right\},\\
\displaystyle&q=(\Delta-\gamma\mu^{2})\psi+\frac{2\mu}{\mathrm{Ro}}-\mu h_{1}+\mu b,\\
\displaystyle&j=h_{1}-\mu\psi.
\end{aligned}
\right.
\end{equation}

\section{Hamiltonian formulation of TQG}
\label{sect3}
In this section, we give a Hamiltonian formulation of the spherical TQG model \eqref{TQGsphere1} in terms of a non-canonical Lie--Poisson semidirect product bracket. Indeed, system \eqref{TQGsphere1} can be seen as a Lie--Poisson flow on the dual of the infinite-dimensional semidirect product Lie algebra $\mathfrak{f}=\mathfrak{f}_{1}\ltimes\mathfrak{f}_{2}$, where $\mathfrak{f}_{1}$ is the space for the potential vorticity field $q$, and $\mathfrak{f}_{2}$ is the space for the buoyancy $b$.

First, we observe that equations \eqref{TQGsphere1} closely resemble the reduced MHD equations \cite{Strauss1976}. Particularly, the potential vorticity $q$ plays the same role as the plasma vorticity, and the buoyancy field $b$ plays the role of the magnetic potential. Furthermore, the thermal (buoyancy) effects on the potential vorticity transport in \eqref{TQGsphere1} are the same as those of the Lorentz force on the vorticity advection in reduced MHD. Additionally, the buoyancy field is transported by the fluid stream function $\psi$, as is the magnetic potential in reduced MHD.

One can write the Hamiltonian functional for \eqref{TQGsphere1} as follows:
\begin{equation}
\begin{split}
\label{hamTQG}
H&=\frac{1}{2}\int\limits_{S^{2}}\left(q-\frac{2\mu}{\mathrm{Ro}}+\mu(h_{1}-b)\right)\left(\Delta-\gamma\mu^{2}\right)^{-1}\left(q-\frac{2\mu}{\mathrm{Ro}}+\mu(h_{1}-b)\right)\mathrm{d}x\mathrm{d}y+\int\limits_{S^{2}}bh_{1}\mathrm{d}x\mathrm{d}y={}\\&=\frac{1}{2}\int\limits_{S^{2}}\left(q-\frac{2\mu}{\mathrm{Ro}}+\mu(h_{1}-b)\right)\psi\mathrm{d}x\mathrm{d}y+\int\limits_{S^{2}}bh_{1}\mathrm{d}x\mathrm{d}y.
\end{split}
\end{equation}
The \textit{semidirect product Lie--Poisson bracket} on $\mathfrak{f}^{*}$ was originally introduced for reduced magnetohydrodynamics \cite{MoGr1980,HazMorr,HolmKuper1983, ripa1995low}, and for functionals $F,G$ it reads as follows:
\begin{equation}
\label{TQGbracket}
\begin{split}
 \llbracket F,G\rrbracket&=\int\limits_{S^{2}}\left[q\left\{\frac{\delta F}{\delta q},\frac{\delta G}{\delta q}\right\}+
 b\left(\left\{\frac{\delta F}{\delta b},\frac{\delta G}{\delta q}\right\}+\left\{\frac{\delta F}{\delta q},\frac{\delta G}{\delta b}\right\}\right)\right]\mathrm{d}x\mathrm{d}y.
 \end{split}
\end{equation}
It allows to express the system \eqref{TQGsphere} as
\begin{equation*}
\dot F=\llbracket H,F\rrbracket,
\end{equation*}
where $F$ is a functional depending on $q$ and $b$.

The Hamiltonian formulation of \eqref{TQGsphere1} reveals the following infinite collection of conservation laws, given by
\begin{equation}
\label{TQGcas}
\mathcal{C}_{f}=\int\limits_{S^{2}}f(b)\mathrm{d}x\mathrm{d}y,\quad\mathcal{I}_{g}=\int\limits_{S^{2}}q g(b)\mathrm{d}x\mathrm{d}y,
\end{equation}
where $f$ and $g$ are arbitrary smooth functions. The quantities \eqref{TQGcas} are \textit{Casimirs}, i.e. they commute with any functional $\mathcal{J}(q,b)$ in the sense of the semidirect product bracket \eqref{TQGbracket}: $\llbracket \mathcal{C}_{f},\mathcal{J}\rrbracket=\llbracket \mathcal{I}_{g},\mathcal{J}\rrbracket=0$. Again, the analogy with the reduced MHD models is evident: the Casimir $\mathcal{C}_{f}$ is the \textit{magnetic helicity} in MHD and reflects that the field $b$ is transported; the Casimir $\mathcal{I}_{g}$ is the counterpart of \textit{cross-helicity} in MHD.

\section{Lie--Poisson discretization of TQG on the sphere}
\label{sect4}
The formulation of the spherical TQG model in terms of non-canonical Hamiltonian structures presented in the previous section suggests that a discretized model should ideally preserve fundamental conservation laws such as Casimirs and energy (Hamiltonian).
A numerical method that preserves the Lie--Poisson structure of the governing equations is referred to as a \textit{Lie--Poisson integrator}.
Such structure-preserving integrators typically boast good stability properties, do not suffer from artificial dissipation, and produce physically relevant results by ensuring that the numerical solution is compatible with known physical and mathematical principles. 

The spatial discretization used in the present study is based on replacing the infinite-dimensional Lie--Poisson structure by a finite-dimensional Lie--Poisson structure. This method is due to Zeitlin \cite{Ze1991,Zeit,Ze2005}, who developed a self-consistent finite mode truncation of the ideal Euler equations on the flat torus referred to as \textit{Euler--Zeitlin equations}. The key idea underlying his approach is to approximate the infinite-dimensional Poisson algebra of smooth functions $(C^{\infty}(S^{2}),\left\{\cdot,\cdot\right\})$ by a sequence of Lie algebras of skew-hermitian matrices $\mathfrak{u}(N)$ converging to $(C^{\infty}(S^{2}),\left\{\cdot,\cdot\right\})$ in a certain sense as $N\to+\infty$, see \cite{ChPolt2017}. The Lie algebra structure is provided by the scaled matrix commutator $[\cdot,\cdot]_{N}=
\frac{1}{\hbar}[\cdot,\cdot]$, where $\hbar=2/\sqrt{N^{2}-1}$.  

A fully structure-preserving discretization is subsequently achieved by choosing an appropriate time integration method. Previous studies have adopted an isospectral integrator \cite{ModViv} to accompany Zeitlin's model for the 2D Euler equations on the sphere. The resulting discrete model fully preserves the underlying Lie--Poisson geometry of spherical hydrodynamics and has been utilized for the studies of ideal turbulence \cite{ModViv1,ModViv2,CiViMo2023}, as well as of QG turbulence \cite{LFEG2024, franken2024singlelayer, franken2025casimir}.

A crucial difference between the Euler and QG models compared to the TQG model is that the dynamics of the latter involves multiple fields. The presence of the buoyancy contributes to the potential vorticity evolution such that it is no longer a single transport equation. Correspondingly, the enstrophy is not preserved in TQG, as opposed to the Euler and QG equations. An analogous result holds for the Zeitlin approximation of the TQG model. Namely, the matrix flow is isospectral only for the buoyancy matrix, while the spectrum of the potential vorticity matrix is not preserved. This means that the isospectral integrator of \cite{ModViv} needs to be extended to be compatible with semidirect product Lie algebras, which was achieved in \cite{ModRoop}. Here, we apply this method to simulate the TQG model on the sphere. 

The finite-dimensional Lie--Poisson structure is obtained via a projection $p_{N}\colon C^{\infty}(S^{2})\to\mathfrak{u}(N)$ from smooth functions on the sphere to skew-hermitian matrices. The projection operator acts on smooth functions (e.g., the potential vorticity $q\in C^{\infty}(S^{2})$) via their spherical harmonic expansion $q=\sum_{l=0}^{\infty}\sum_{m=-l}^{l}q_{lm}Y_{lm}$. Namely, truncating this decomposition at some degree $N-1$ and replacing the spherical harmonics $Y_{lm}$ with \textit{matrix harmonics} $T_{lm}(N)$, we obtain the \textit{potential vorticity matrix}
\begin{equation*}
Q=p_{N}(q)=\sum\limits_{l=0}^{N-1}\sum\limits_{m=-l}^{l}q_{lm}(\mathrm{i}T_{lm}).
\end{equation*}
Here $T_{lm}$ are the eigen-matrices of the \textit{Hoppe--Yau Laplacian} $\Delta_{N}$ \cite{HoppYau}, i.e. $\Delta_{N}T_{lm}=-l(l+1)T_{lm}$, and $\mathrm{i}$ is the imaginary unit. 
Conversely, one can reconstruct $q$ up to degree $N$ using the spectral coefficients of $Q$. The stream function and buoyancy are similarly approximated by a stream matrix and buoyancy matrix.
The projection operator is applied extensively to products of functions throughout the numerical method. Given two functions $f$ and $g$ with projections $F=p_{N}(f)$ and $G=p_{N}(g)$, the projection of the product $fg$ is given by
\begin{equation*}
p_{N}(fg)=-\frac{\mathrm{i}}{2}\sqrt{\frac{N}{4\pi}}(FG+GF),
\end{equation*}
which converges to $fg$ at a rate $\mathcal{O}(1/N)$. Henceforth, we denote the projection of the product by $F\odot G=p_N(fg)$. For a detailed analysis, we refer to \cite{franken2024singlelayer}.

We define $Q=p_N(q), P=p_N(\psi), B=p_N(b), M=p_N(\mu), S=p_N(\mu^2)$ and $ H_1=p_N(h_1)$. Then the finite-dimensional matrix approximation of \eqref{TQGsphere1} reads
\begin{equation}
\label{TQGquant}
    \left\{\begin{aligned}
        &\dot Q=[Q, P]_{N} + [B,J]_{N} \\
        &\dot B=[B, P]_{N} \\
        &\Delta_NP - \gamma S\odot P= Q - M\odot(B-H_1)-\frac{2M}{\mathrm{Ro}} \\
        &J=H_1-M\odot P
    \end{aligned}
    \right.
\end{equation}

System \eqref{TQGquant} will be referred to as the \textit{TQG--Zeitlin} equations, and forms a finite-dimensional Lie--Poisson system on the dual $\mathfrak{f}^{*}$ of the semidirect product Lie algebra $\mathfrak{f}=\mathfrak{u}(N)\ltimes\mathfrak{u}(N)^{*}$. The Hamiltonian is given by
\begin{equation}
\label{qHamiltTQG}
H(Q,B)=\frac{1}{2}\mathrm{tr}\left(\left(Q-\frac{2M}{\mathrm{Ro}}+M\odot(H_{1}-B)\right)^{\dag}P\right)+\mathrm{tr}(B^{\dag}H_{1}),
\end{equation}
and the Casimirs are
\begin{equation}
\label{qCasTQG}
\mathcal{C}_{f}^{N}=\frac{4\pi}{N}\mathrm{tr}\left(f(B)\right),\quad \mathcal{I}_{g}^{N}=\frac{4\pi}{N}\mathrm{tr}(Qg(B)).
\end{equation}
The Hamiltonian \eqref{qHamiltTQG} and the Casimirs \eqref{qCasTQG} with the functions $f$ and $g$ chosen to be monomials, converge to their continuous analogs \eqref{hamTQG} and \eqref{TQGcas} as $N\to\infty$, see \cite{ModRoop}.

Finally, a suitable structure-preserving Lie--Poisson time integrator for the flow \eqref{TQGquant} has been developed in \cite{ModRoop} and is adopted in the present study. It preserves the Casimirs \eqref{qCasTQG} for the TQG--Zeitlin equations  up to machine precision and nearly preserves the energy \eqref{hamTQG}. For a time step size $h$, the scheme $\Phi_{h}(Q_{n},B_{n})\mapsto(Q_{n+1},B_{n+1})$ is given by
\begin{equation}
\label{TQGmethod1}
\begin{aligned}
&B_{n}=\tilde B-\frac{h}{2}[\tilde B,\tilde{P}]-\frac{h^{2}}{4}\tilde{P}\tilde B\tilde{P}, \\
&B_{n+1}=B_{n}+h[\tilde B,\tilde{P}],\\
&Q_{n}=\tilde Q-\frac{h}{2}[\tilde Q,\tilde{P}]-\frac{h}{2}[\tilde B,\tilde J]-\frac{h^{2}}{4}\left(\tilde{P}\tilde Q\tilde{P}+\tilde J\tilde B\tilde{P}+\tilde{P}\tilde B\tilde J\right), \\
&Q_{n+1}=Q_{n}+h[\tilde Q,\tilde P]+h[\tilde B,\tilde J],
\end{aligned}
\end{equation}
where $\tilde{P}$ is computed from $\tilde{Q}$ via the inhomogeneous Helmholtz equation.
The subscripts $n$, $n+1$ denote the time instances. 
The method \eqref{TQGmethod1} is a midpoint method that preserves the Lie--Poisson structure and has order of consistency $\mathcal{O}(h^{2})$. An implicit system of nonlinear equations is solved via fixed point iteration to find the variables with a tilde $\tilde{\cdot}$ at the midpoint of each time step.
An important advantage of the scheme \eqref{TQGmethod1} is the absence of group-to-algebra maps typical for Lie--Poisson integrators, which makes the method efficiently applicable to high-dimensional Lie--Poisson matrix flows.

The computational cost of the presented algorithm is dominated by the computation of the commutator, which involves dense matrix multiplications and requires $\mathcal{O}(N^3)$ operations. The inhomogeneous Helmholtz problem for computing the stream matrix $P$ is solved efficiently using a diagonal splitting approach as described by \cite{franken2024singlelayer}, in $\mathcal{O}(N^2)$ operations.
\section{Simulation results}
\label{sect5}

We proceed to illustrate several examples of structure-preserving simulation of the TQG model on the sphere, where we concentrate on the conservation properties of the numerical method and qualitative flow dynamics. A distinctive feature of the Lie--Poisson integrator described in the previous section is the absence of numerical (artificial) viscosity. This enables indefinitely long simulations even without external forcing and viscous dissipation or other types of regularization. Furthermore, the absence of external forcing and dissipation ensures that the Casimirs are preserved and the energy is nearly preserved.

The presented numerical results are used to demonstrate the capabilities of the computational method and to establish that simulation of a large-scale TQG model on the sphere is achievable. We adopt a Rossby number $\mathrm{Ro}=0.01$, which is relevant for planetary motions of the Earth's atmosphere (see, e.g., \cite{LFEG2024}). The Lamb parameter takes the values $\gamma=100$ and $\gamma=1000$, which correspond to different rotation frequencies. These parameters are chosen to illustrate qualitative flow features, rather than exactly match Earth's dynamics.  Explicit comparison with Earth's dynamics requires external forcing and dissipation and warrants a separate study. 

\subsection{Trivial bottom topography}
The first set of simulations is performed with a trivial bottom topography, $h_{1}=0$. The Rossby number is $\mathrm{Ro}=0.01$, and the Lamb parameter $\gamma=100$. The matrix dimension is $N=512$, and the final simulation time is $T=33$. Initial distributions for the potential vorticity $q$ and the buoyancy $b$ are randomly generated smooth fields. Namely, we generate initial conditions by specifying the first 100 spherical harmonic coefficients $q_{lm}$ of the potential vorticity as an array of samples from the standard normal distribution, while leaving the rest of the coefficients trivial. The initial buoyancy field is generated in the same way.

\begin{figure}[h]
\begin{minipage}[h]{0.5\linewidth}
\center{\includegraphics[scale=0.45]{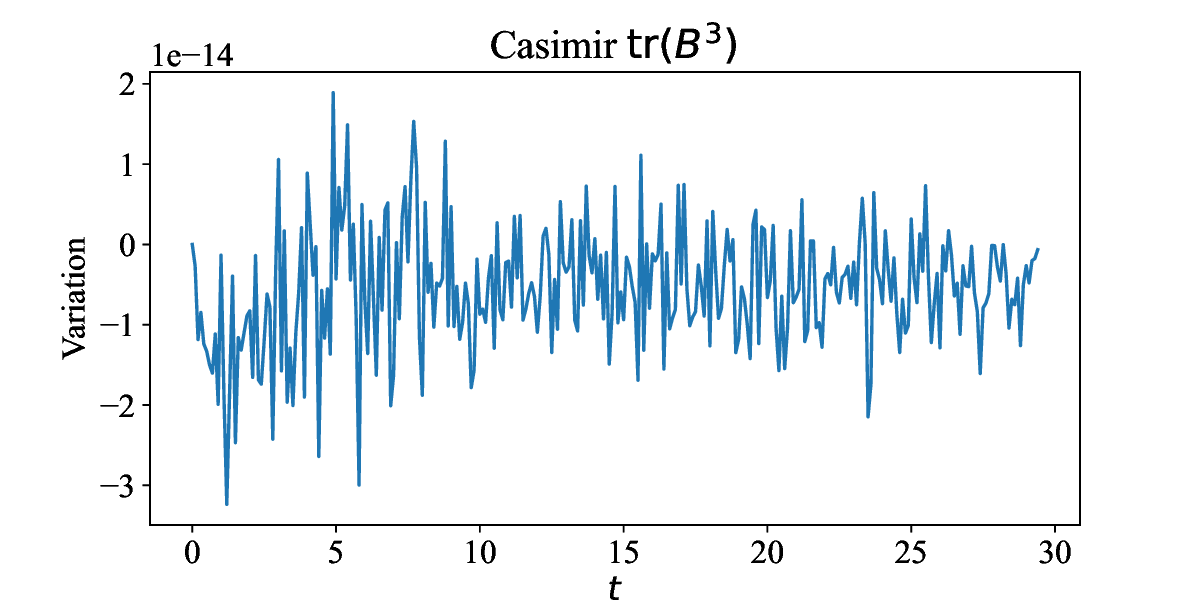}}
\end{minipage}%
\hfill
\begin{minipage}[h]{0.5\linewidth}
\center{\includegraphics[scale=0.45]{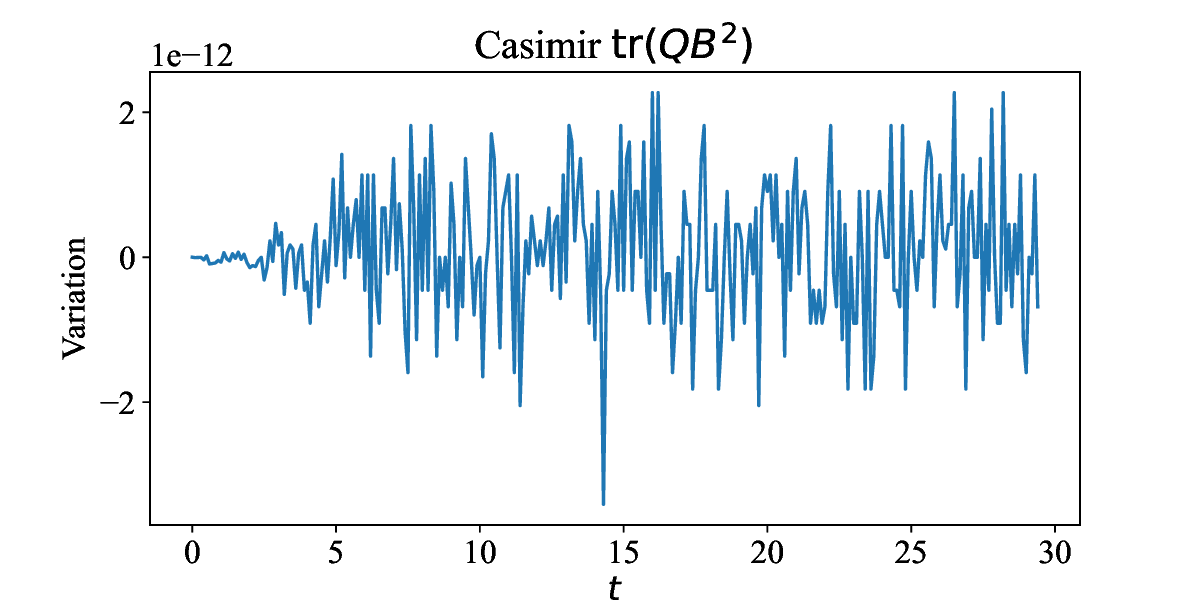}}
\end{minipage}
\hfill
\caption{Variation of Casimirs $\mathrm{tr}(B^{3})$ and $\mathrm{tr}(QB^{2})$, corresponding to the choice $f(B)=B^{3}$ and $g(B)=B^{2}$ in \eqref{qCasTQG}. The magnitude of the variation indicates preservation up to machine precision.}
\label{cas_pres}
\end{figure}

\begin{figure}[ht!]
\centering
\includegraphics[scale=0.5]{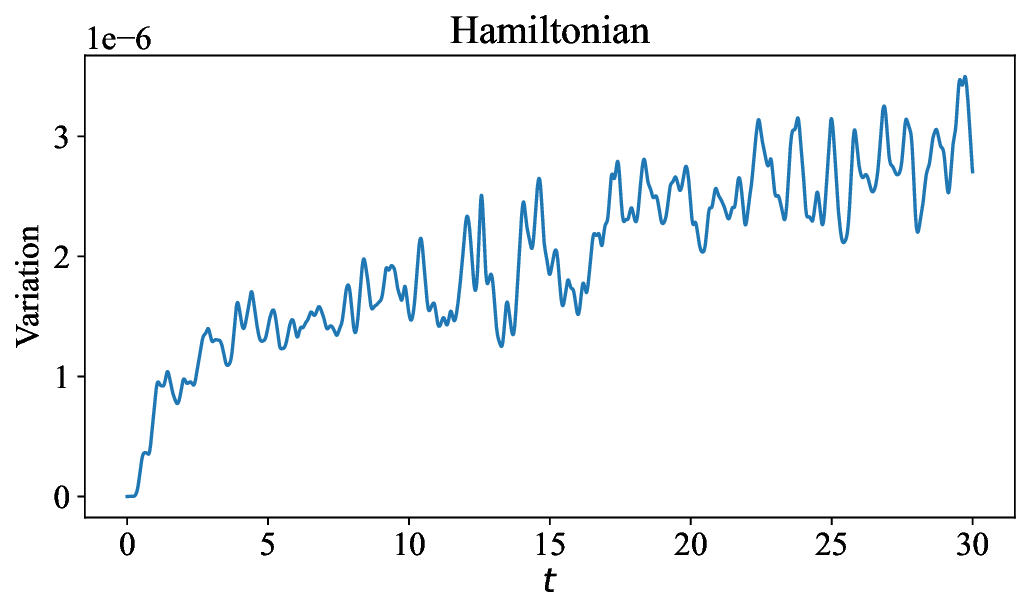}
\caption{Variation of the Hamiltonian \eqref{qHamiltTQG}. The magnitude of the variation indicates near preservation of the Hamiltonian.}
\label{ham_preserv}
\end{figure}

First, we demonstrate that the scheme \eqref{TQGmethod1} exactly preserves the Casimir invariants. In Fig.~\ref{cas_pres}, it is shown that the variation of the Casimirs is of machine precision. Fig.~\ref{ham_preserv} shows near preservation of the Hamiltonian \eqref{qHamiltTQG}, which is validated by the variation magnitude $10^{-6}$. The integrator \eqref{TQGmethod1} ensures exact preservation of Casimirs, but cannot achieve exact preservation of the Hamiltonian function, which is a well-known feature of geometric integrators \cite{HaiLubWan}.

\begin{figure}[ht]
    \begin{minipage}[t]{0.5\textwidth}
    \captionsetup{width=1\textwidth,labelformat=empty} %
      \includegraphics[width=\linewidth]{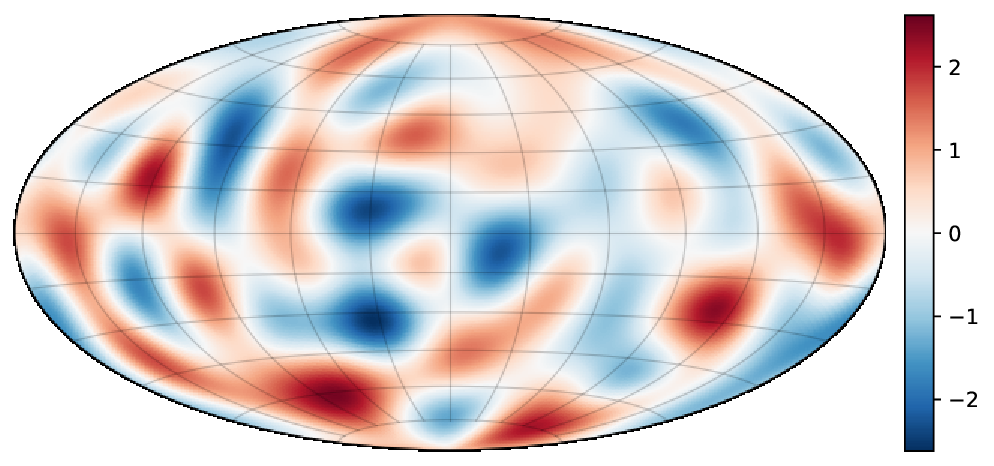} 
      \caption*{Potential vorticity $q$,  $t=0$}
    \end{minipage}%
    \hfill 
    \begin{minipage}[t]{0.5\textwidth}
    \captionsetup{width=1\textwidth,labelformat=empty}
      \includegraphics[width=\linewidth]{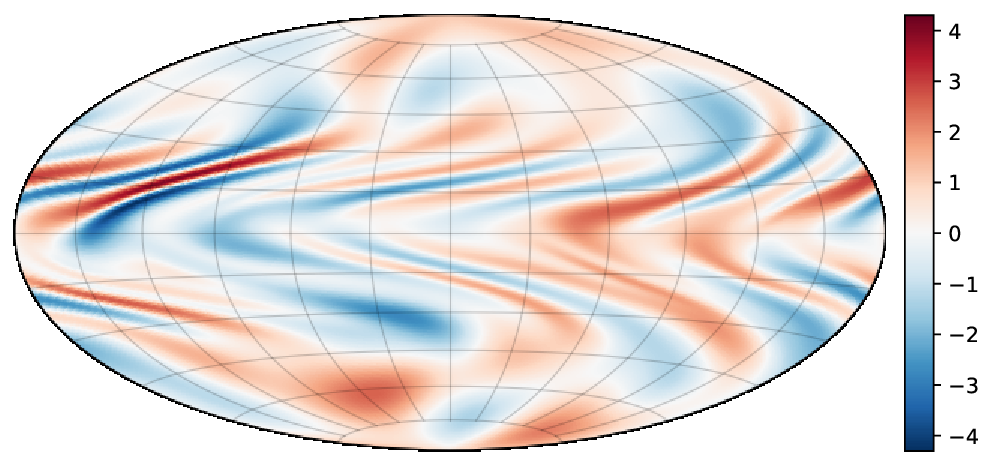} 
      \caption*{Potential vorticity $q$,  $t=1$}
    \end{minipage}%
    \hfill 
    \begin{minipage}[t]{0.5\textwidth}
    \captionsetup{width=1\textwidth,labelformat=empty}
      \includegraphics[width=\linewidth]{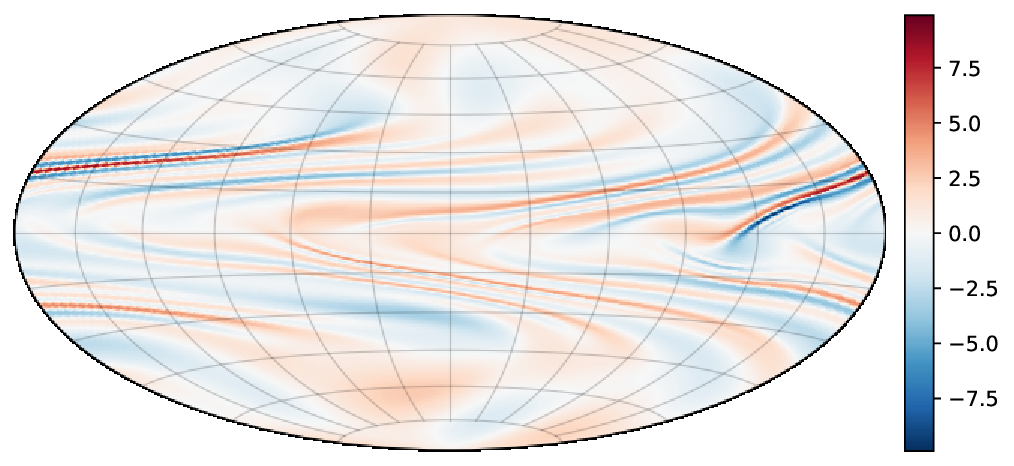} 
      \caption*{Potential vorticity $q$,  $t=2$}
    \end{minipage}%
    \hfill 
    \begin{minipage}[t]{0.5\textwidth}
    \captionsetup{width=1\textwidth,labelformat=empty}
      \includegraphics[width=\linewidth]{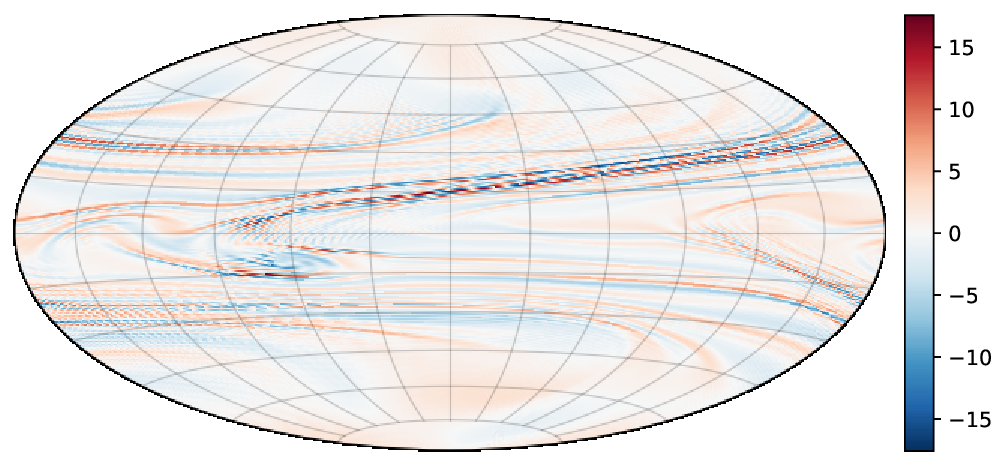} 
      \caption*{Potential vorticity $q$,  $t=4$}
    \end{minipage}
\caption{Evolution of the potential vorticity $q(t)$ field for the spherical TQG equations. Initial smooth randomly generated field develops vorticity filaments and small scale dynamics with simultaneous growth of the magnitude.}
\label{q_field}
\end{figure}

The simulation results for the potential vorticity $q$ are shown in Fig.~\ref{q_field}. The potential vorticity is shown on shorter time scales than the buoyancy. The presence of two Poisson brackets in the potential vorticity equation causes growth and rapid development of small-scale features, after which it becomes difficult to discern flow patterns. That is, buoyancy causes amplification of the potential vorticity, visible in the amplitude of $q$. We provide the dynamics up to $t=4$, when the dynamics are still well-resolved for the chosen spatial resolution. The effect of rotation is clearly visible near the equator, where elongated vorticity filaments with sharp gradients form.

\begin{figure}[ht]
    \begin{minipage}[t]{0.5\textwidth}
    \captionsetup{width=1\textwidth,labelformat=empty} %
      \includegraphics[width=\linewidth]{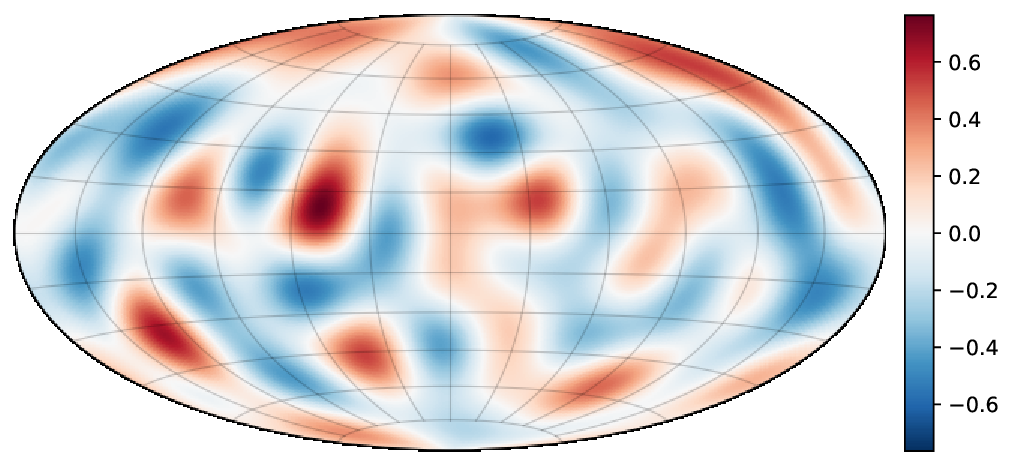} 
      \caption*{Buoyancy $b$,  $t=0$}
    \end{minipage}%
    \hfill 
    \begin{minipage}[t]{0.5\textwidth}
    \captionsetup{width=1\textwidth,labelformat=empty}
      \includegraphics[width=\linewidth]{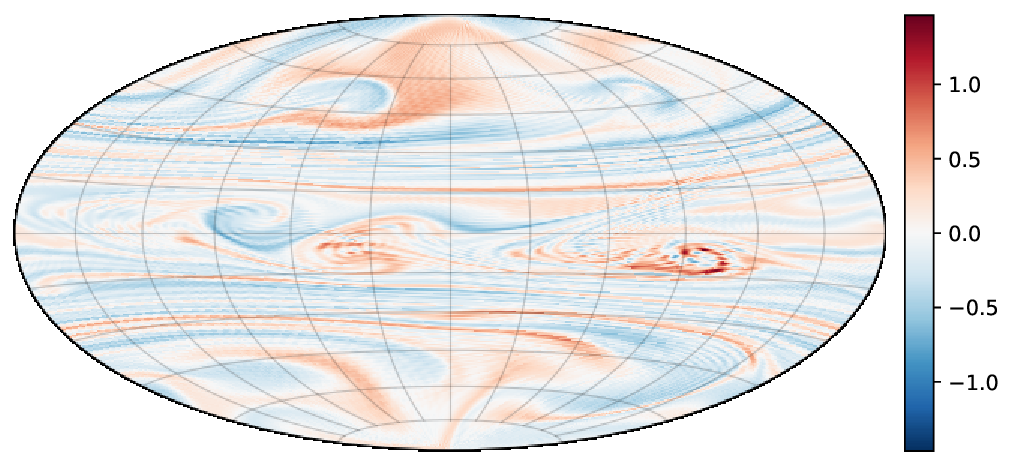} 
      \caption*{Buoyancy $b$,  $t=7$}
    \end{minipage}%
    \hfill 
    \begin{minipage}[t]{0.5\textwidth}
    \captionsetup{width=1\textwidth,labelformat=empty}
      \includegraphics[width=\linewidth]{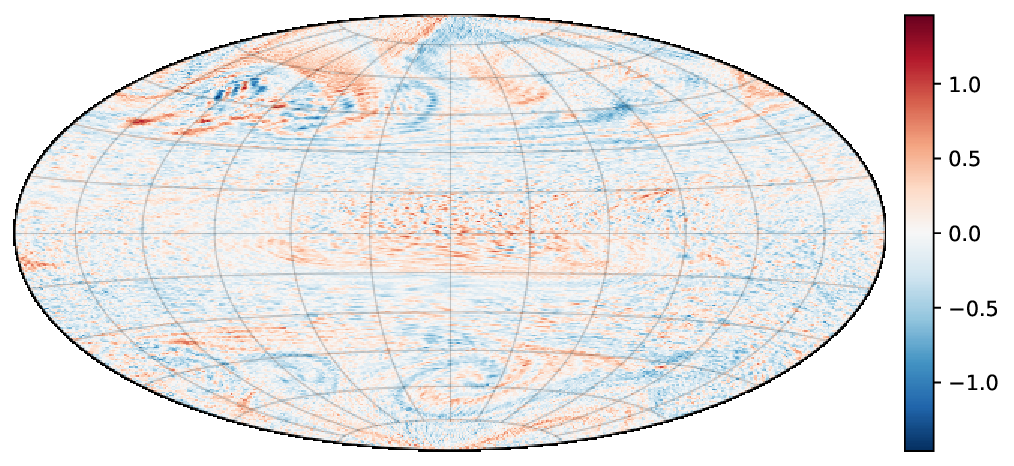} 
      \caption*{Buoyancy $b$,  $t=18$}
    \end{minipage}%
    \hfill 
    \begin{minipage}[t]{0.5\textwidth}
    \captionsetup{width=1\textwidth,labelformat=empty}
      \includegraphics[width=\linewidth]{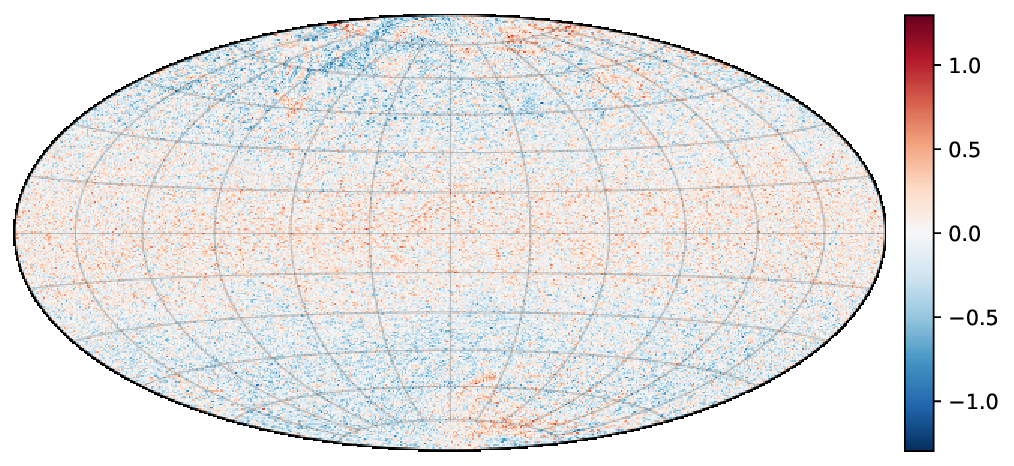} 
      \caption*{Buoyancy $b$,  $t=33$}
    \end{minipage}
\caption{Evolution of the buoyancy $b(t)$ field for the spherical TQG equations. Initially smooth field develops turbulent mixing in the equatorial domain along with zonal structures in mid-latitudes. Final distribution has a zonal structure with generally positive buoyancy at the equator and generally negative at the poles.}
\label{b_field}
\end{figure}

The evolution of the buoyancy is illustrated in Fig.~\ref{b_field}. Turbulent mixing is observed closest to the equator, where the Coriolis effects are small. Simultaneously, in mid-latitudinal bands, elongated zonal structures develop under the effect of rotation.
Gradually, small-scale features develop in the mid-latitudes as well within larger zonal buoyancy structures. The zonal structures ultimately coalesce to span across the equator in the long-time distribution of the buoyancy, showing generally positive values at the equator and negative values at the poles. 

The formation of small scales ultimately leads to a `noisy' solution, observed in, e.g., the last depicted snapshot of the buoyancy. This behavior is native to the Zeitlin approximation and accords with earlier observations of the potential vorticity in the freely evolving two-dimensional Euler equations \cite{ModViv1} and QG equations \cite{LFEG2024}. Specifically, this is a consequence of the conservation properties of the numerical method. Nonlinear advection ensures distribution of energy over all resolvable scales of motion. Eventually, the energy distribution stabilizes and a statistically steady flow state is reached. The absence of any viscous or numerical (artificial) dissipation ensures that the small-scale features persist, even after arbitrarily long simulation times. Viscous dissipation suppresses the small-scale energy in the Zeitlin approximation as was previously demonstrated for the two-dimensional Navier--Stokes equations \cite{cifani2022casimir} and geostrophic turbulence \cite{franken2024singlelayer}, respectively leading to recognizable flow structures associated with isotropic turbulence and zonal jet formation. Numerical simulations of viscous Zeitlin-TQG turbulence are the subject ongoing study.

\begin{figure}[h]
\begin{minipage}[h]{0.5\linewidth}
\center{\includegraphics[scale=0.45]{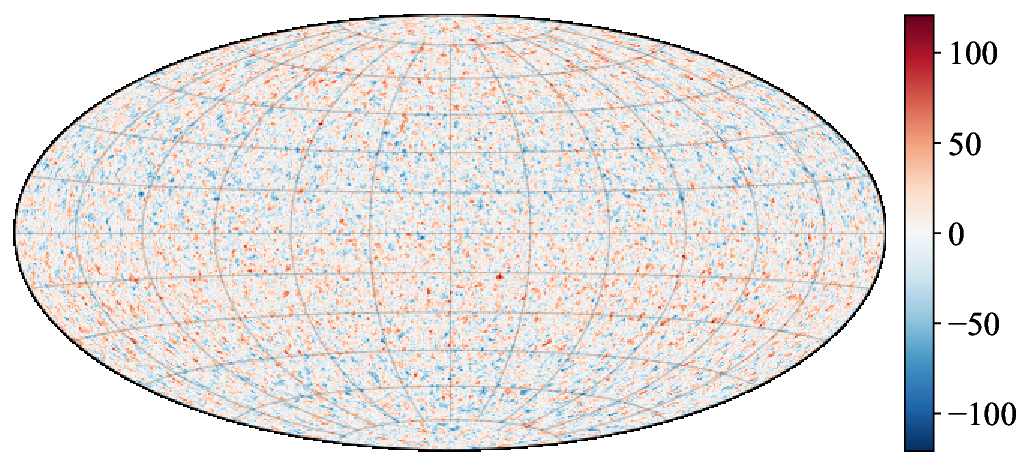}}
\end{minipage}%
\hfill
\begin{minipage}[h]{0.5\linewidth}
\center{\includegraphics[scale=0.45]{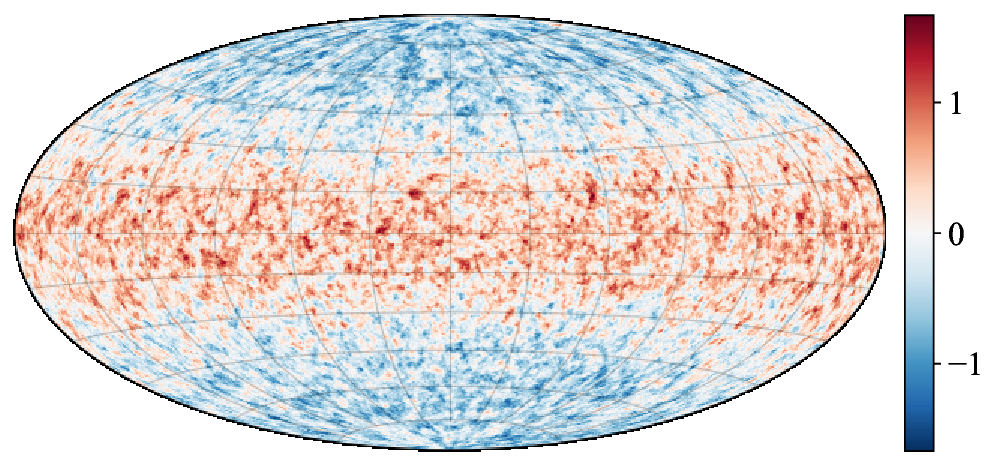}}
\end{minipage}
\hfill
\caption{Final distributions of the potential vorticity (left) and buoyancy (right) after applying the Helmholtz filter.}
\label{qb_filter}
\end{figure}

In Fig.~\ref{qb_filter}, we show the final distributions of the potential vorticity and the buoyancy. In these figures, we apply a low-pass Helmholtz filter $(1-\alpha^2\Delta)$ with $\alpha=1/64$ to illustrate large-scale zonal flow features present in the long-time solution. The filter is solely for visualization purposes and is not used in the numerical simulation.

\subsection{Buoyancy-bathymetry interaction} 
\begin{figure}[h]
    \begin{minipage}[t]{0.33\textwidth}
    \captionsetup{width=1\textwidth,labelformat=empty} %
      \includegraphics[width=\linewidth]{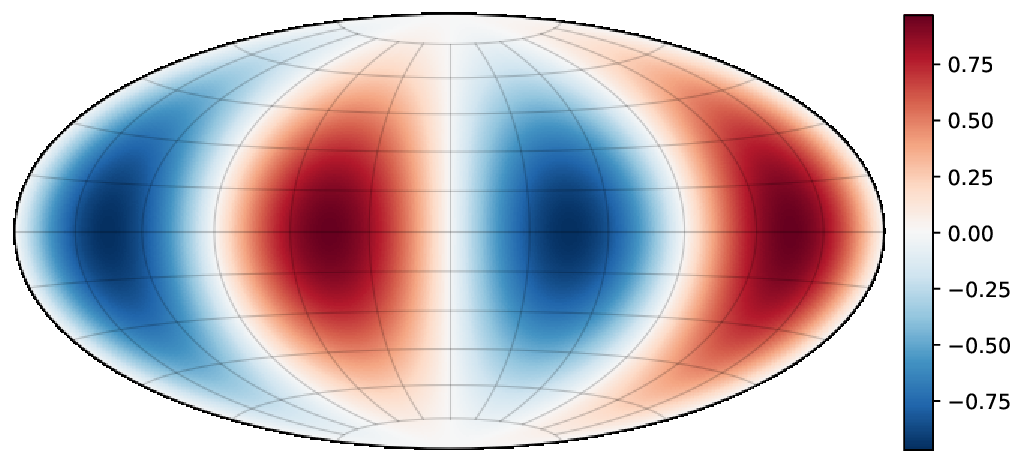} 
    \end{minipage}%
    \hfill 
    \begin{minipage}[t]{0.33\textwidth}
    \captionsetup{width=1\textwidth,labelformat=empty}
      \includegraphics[width=\linewidth]{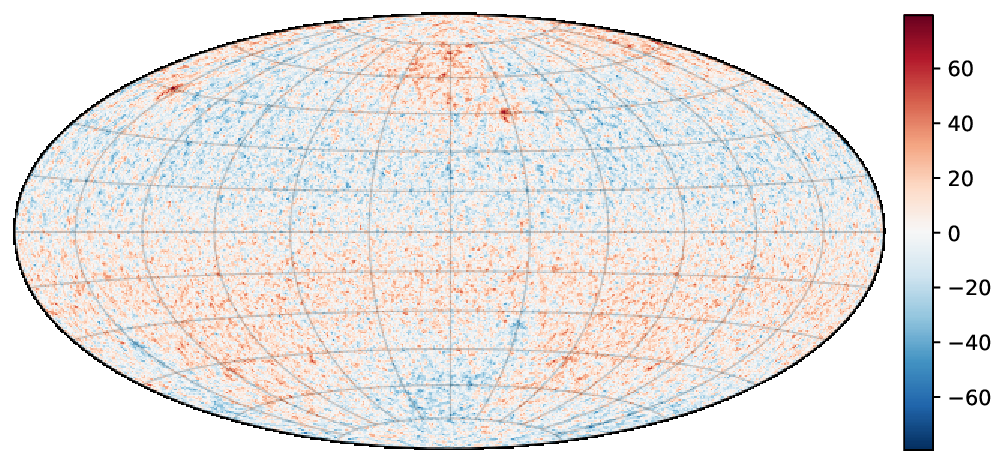} 
    \end{minipage}%
    \hfill 
    \begin{minipage}[t]{0.33\textwidth}
    \captionsetup{width=1\textwidth,labelformat=empty}
      \includegraphics[width=\linewidth]{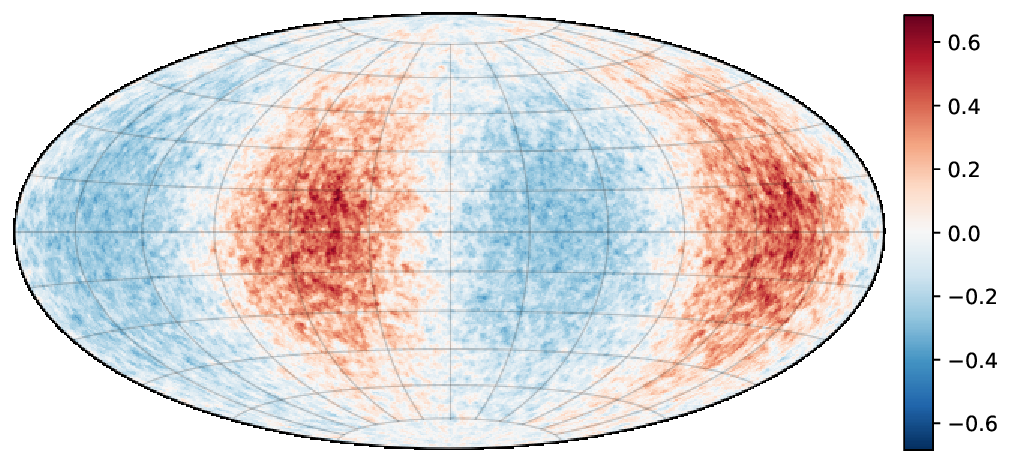} 
    \end{minipage}%
\caption{Adopted bathymetry profile (left) and final distributions of the potential vorticity (center) and buoyancy (right). The potential vorticity exhibits zonal structure and aligns with the gradient of the bathymetry, while the buoyancy aligns with the bathymetry. }
\label{long_time_secondcase}
\end{figure}

\begin{figure}[ht]
    \begin{minipage}[t]{0.5\textwidth}
    \captionsetup{width=1\textwidth,labelformat=empty} %
      \includegraphics[width=\linewidth]{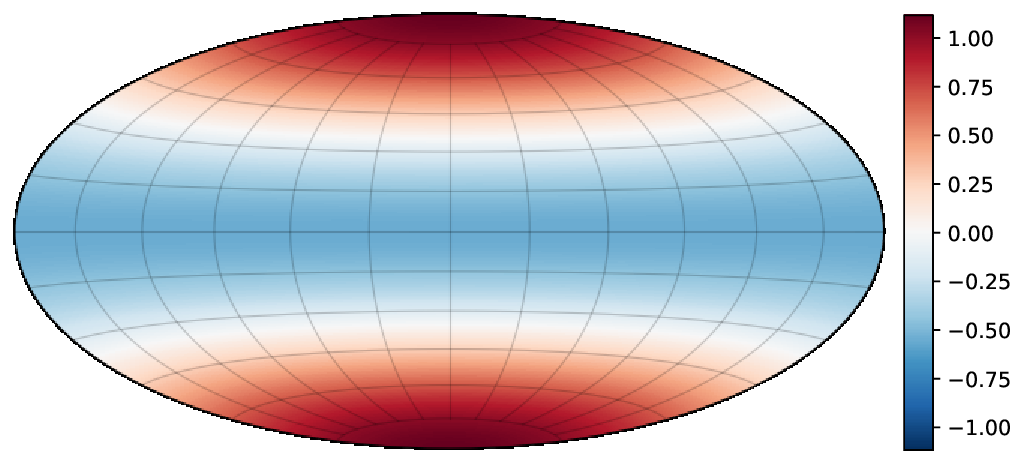} 
      \caption*{Buoyancy $b$,  $t=0.1$}
    \end{minipage}%
    \hfill 
    \begin{minipage}[t]{0.5\textwidth}
    \captionsetup{width=1\textwidth,labelformat=empty}
      \includegraphics[width=\linewidth]{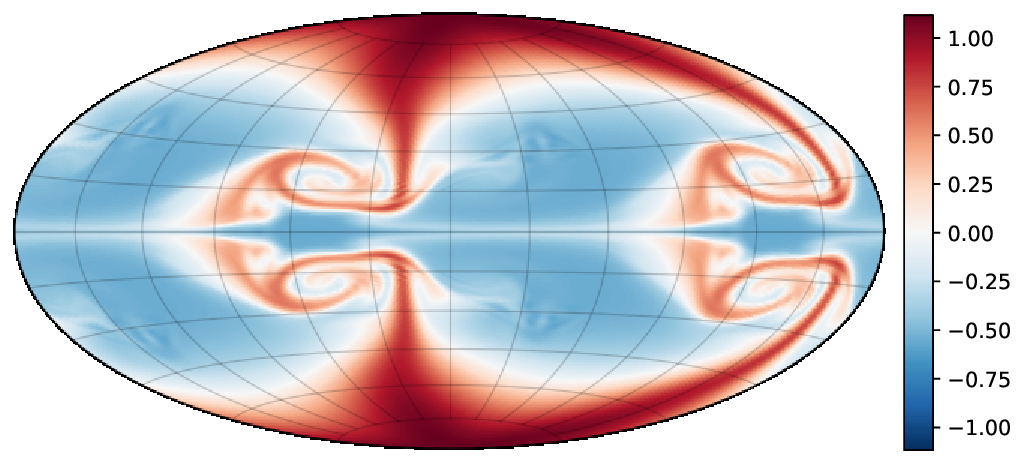} 
      \caption*{Buoyancy $b$,  $t=8.4$}
    \end{minipage}%
    \hfill 
    \begin{minipage}[t]{0.5\textwidth}
    \captionsetup{width=1\textwidth,labelformat=empty}
      \includegraphics[width=\linewidth]{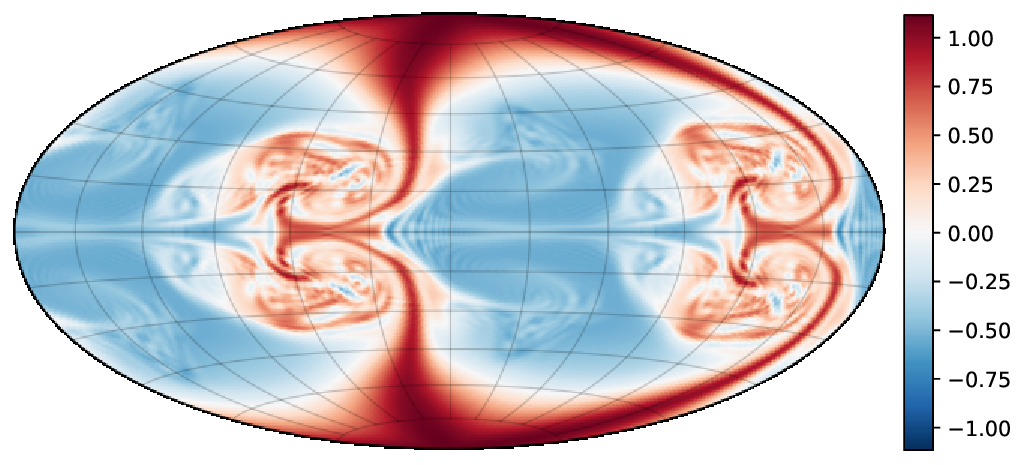} 
      \caption*{Buoyancy $b$,  $t=10$}
    \end{minipage}%
    \hfill 
    \begin{minipage}[t]{0.5\textwidth}
    \captionsetup{width=1\textwidth,labelformat=empty}
      \includegraphics[width=\linewidth]{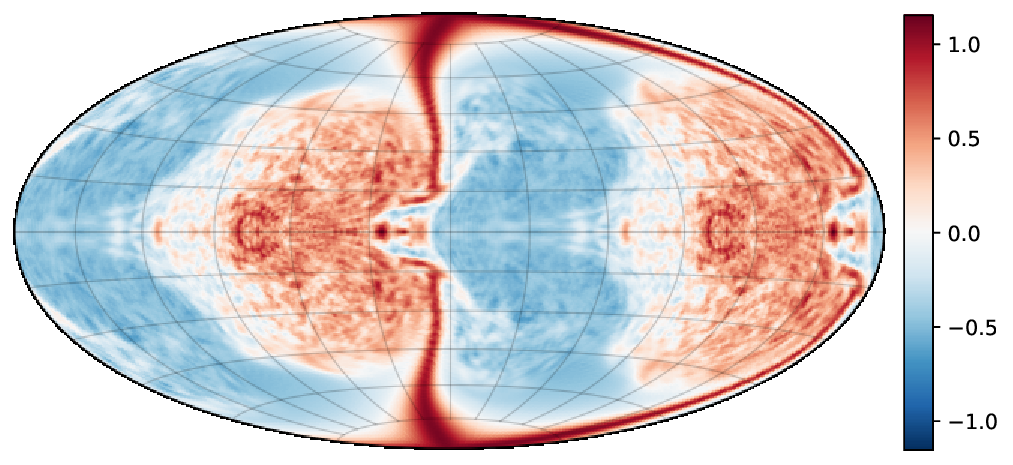} 
      \caption*{Buoyancy $b$,  $t=14$}
    \end{minipage}
\caption{Buoyancy dynamics. From an initially smooth field, mushroom-like dipole structures form and roll up asymmetrically under the influence of rotation and the interaction with the bathymetry. Circulation is induced via the buoyancy-bathymetry interaction and turbulence is generated which eventually aligns with the bathymetry.}
\label{b_field_secondcase}
\end{figure}

The second test case is designed to exhibit how circulation is created by the misalignment of the gradients of the buoyancy and the bathymetry. The bathymetry and buoyancy are initialized as smooth fields, respectively depicted in the left panel of Fig.~\ref{long_time_secondcase} and the top left panel of Fig.~\ref{b_field_secondcase}. The initial vorticity is set to zero. We adopt a Rossby number $\mathrm{Ro}=0.01$ and Lamb parameter $\gamma=1000$.

Previous numerical results of TQG in the $\beta$-plane approximation \cite{HolmLuesinkPan2021} suggest emergent sequences of coherent flow patterns at successively smaller scales at the onset of high-wavenumber instabilities. Consequently, scale-resolving numerical simulations quickly become challenging without regularization and necessitate, e.g., viscous dissipation to specify a smallest length scale of the flow. However, numerical simulations of relevant hyperbolic equations such as the two-dimensional Euler or QG equations eventually face the same challenge, but can instead still be used to study large-scale flow patterns (see \cite{ModViv2, LFEG2024, franken2025casimir}). In the present study, the hyperbolic TQG equations are simulated without regularization and hence the evolution of small scales is impeded by numerical resolution and leads to `noisy' solutions. As described before, we therefore only illustrate large-scale flow features by applying a low-pass Helmholtz filter as a post-processing step. 

We first highlight the evident influence of the bathymetry on the long-time distributions of the potential vorticity and buoyancy, shown in Fig.~\ref{long_time_secondcase}. The bathymetry profile is clearly observed in the large-scale buoyancy. This implies an alignment of the gradients of the respective fields, minimizing the circulation induced by $\{b, j\}$ in the potential vorticity dynamics. A zonal structure is observed in the potential vorticity, similar to the first test case. However, an imprint of the bathymetry is visible, suggesting an alignment of the potential vorticity with the \textit{gradient} of the bathymetry in regions where the zonal velocity is not dominant.

\begin{figure}[ht]
    \begin{minipage}[t]{0.5\textwidth}
    \captionsetup{width=1\textwidth,labelformat=empty} %
      \includegraphics[width=\linewidth]{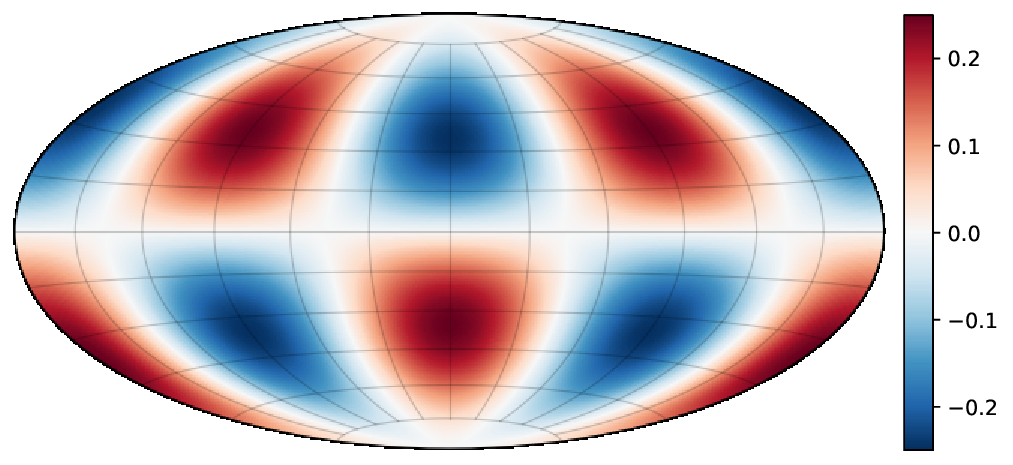} 
      \caption*{Potential vorticity $q$,  $t=0.1$}
    \end{minipage}%
    \hfill 
    \begin{minipage}[t]{0.5\textwidth}
    \captionsetup{width=1\textwidth,labelformat=empty}
      \includegraphics[width=\linewidth]{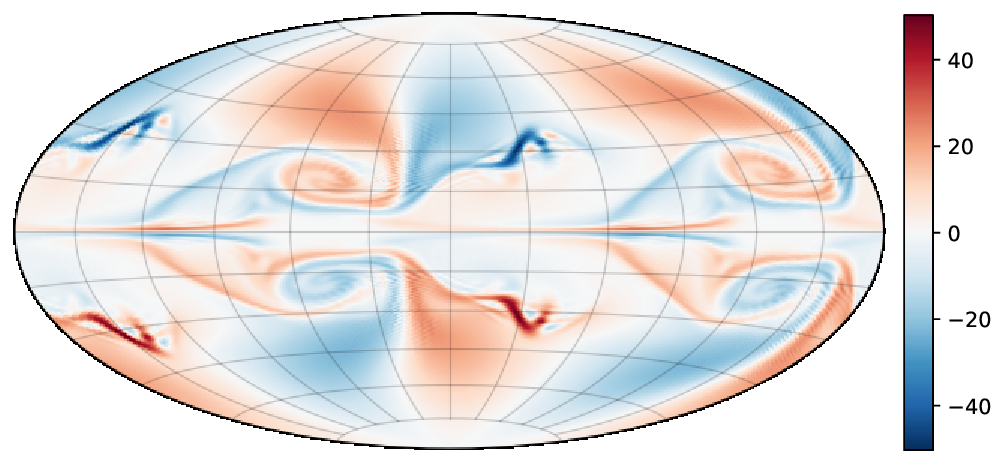} 
      \caption*{Potential vorticity $q$,  $t=8.4$}
    \end{minipage}%
    \hfill 
    \begin{minipage}[t]{0.5\textwidth}
    \captionsetup{width=1\textwidth,labelformat=empty}
      \includegraphics[width=\linewidth]{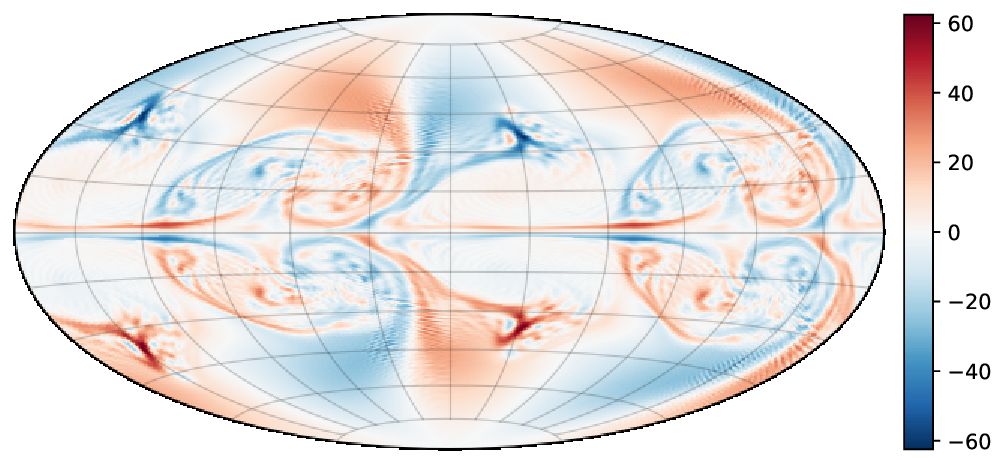} 
      \caption*{Potential vorticity $q$, $t=10$}
    \end{minipage}%
    \hfill 
    \begin{minipage}[t]{0.5\textwidth}
    \captionsetup{width=1\textwidth,labelformat=empty}
      \includegraphics[width=\linewidth]{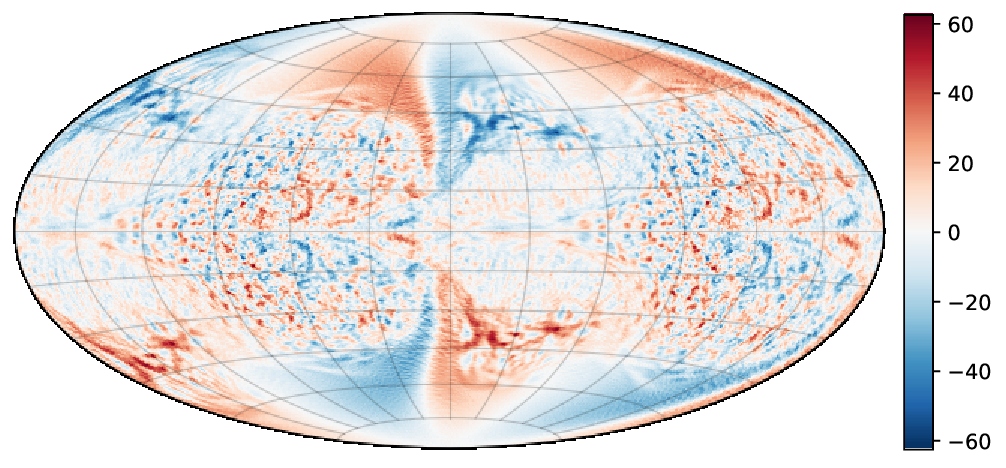} 
      \caption*{Potential vorticity $q$,  $t=14$}
    \end{minipage}
\caption{Potential vorticity dynamics. The misalignment of buoyancy and bathymetry induces circulation in an initially trivial potential vorticity field, which subsequently develops dipole-like circulations and turbulent flow near the equator.}
\label{q_field_secondcase}
\end{figure}

The evolution of the potential vorticity and buoyancy are shown in Fig.~\ref{q_field_secondcase} and Fig.~\ref{b_field_secondcase}, respectively, depicting the filtered fields. After a short simulation time, potential vorticity circulation patterns are induced by misalignment of the buoyancy and bathymetry. Additionally, the rotation causes a strong shear in the zonal velocity near the equator. This combination leads to plumes of buoyancy that are advected parallel to the equator. Mushroom-shaped buoyancy dipoles form, similar to the $\beta$-plane case \cite{HolmLuesinkPan2021}, and an asymmetric roll-up occurs due to the latitude-dependent Coriolis force. The plumes form buoyancy fronts while turbulence is generated via the buoyancy-bathymetry interaction.

\section{Conclusion and outlook}
\label{sect6}
In this paper, we have derived the thermal quasi-geostrophic (TQG) model on the sphere, highlighted its Hamiltonian structure and conserved quantities, and provided a structure-preserving computational method for numerical flow simulations. The model derivation followed from asymptotic expansion of the thermal rotating shallow water equations. The full variation of the Coriolis parameter $f$ is retained, with the simplifying assumption $f\nabla\psi\approx\nabla(f\psi)$. This leads to a TQG formulation for the entire sphere, distinct from previously studied $f$-plane and $\beta$-plane approximations. The resulting system shares a semidirect product Lie--Poisson formulation also observed in ideal two-dimensional magnetohydrodynamics (MHD). The identification of the geometric structure readily enabled numerical simulations using a Casimir-preserving Lie--Poisson integrator for semidirect products \cite{ModRoop}. The conservation properties were shown numerically for freely evolving TQG turbulence. Large-scale flow patterns revealed intricate dynamic interplay between potential vorticity, buoyancy, and bathymetry. 

Future research is dedicated to the numerical simulation and assessment of regularized TQG turbulence, to advance the understanding of TQG on the sphere as a model for planetary flow dynamics. A prime example is the inclusion of dissipation, which can serve, e.g., as a sub-mesoscale parametrization. Further work is needed to study the formation and persistence of large-scale coherent structures such as fronts, gyres, and jets. A closer comparison of the presented discretization to non-geometric methods, and a detailed analysis of importance of the geometric discretization is a subject of ongoing research. Comparing spherical TQG dynamics to TQG flow on the $\beta$ plane would be another promising direction of research.

\bibliographystyle{spmpsci}
\bibliography{TQG_main.bib}

\begin{thebibliography}{10}
\providecommand{\url}[1]{{#1}}
\providecommand{\urlprefix}{URL }
\expandafter\ifx\csname urlstyle\endcsname\relax
  \providecommand{\doi}[1]{DOI~\discretionary{}{}{}#1}\else
  \providecommand{\doi}{DOI~\discretionary{}{}{}\begingroup \urlstyle{rm}\Url}\fi

\bibitem{Arn}
Arnold, V.: Sur la g\'eometri\'e diff\'erentielle des groupes de {L}ie de dimension infnie et ses applications \'a l'hydrodynamique des fluides parfaits.
\newblock Ann. Inst. Fourier (Grenoble) \textbf{16}, 319--361 (1966)

\bibitem{beron2021nonlinear}
Beron-Vera, F.: Nonlinear saturation of thermal instabilities.
\newblock Phys. Fluids \textbf{33}(3) (2021)

\bibitem{Blackburn1985}
Blackburn, M.: {Interpretation of Ageostrophic Winds and Implications for Jet Stream Maintenance}.
\newblock J. Atmos. Sci. \textbf{42}, 2604--2620 (1985)

\bibitem{ChPolt2017}
Charles, L., Polterovich, L.: Sharp correspondence principle and quantum measurements.
\newblock St. Petersburg Math. J. \textbf{29}(1), 177--207 (2018)

\bibitem{cifani2022casimir}
Cifani, P., Viviani, M., Luesink, E., Modin, K., Geurts, B.J.: Casimir preserving spectrum of two-dimensional turbulence.
\newblock Phys. Rev. Fluids \textbf{7}(8), L082601 (2022)

\bibitem{CiViMo2023}
Cifani, P., Viviani, M., Modin, K.: {An efficient geometric method for incompressible hydrodynamics on the sphere}.
\newblock J. Comput. Phys. \textbf{473}, 111772 (2023)

\bibitem{Cressman1958}
Cressman, G.: {Barotropic divergence and very long atmospheric waves}.
\newblock Mon. Wea. Rev. \textbf{86}, 293--297 (1958)

\bibitem{CrisanHolmLuesinkMensahPan2023}
Crisan, D., Holm, D., Luesink, E., Mensah, P., Pan, W.: {Theoretical and Computational Analysis of the Thermal Quasi-Geostrophic Model}.
\newblock J. Nonlinear Sci. \textbf{33}, 96 (2023)

\bibitem{Daley1983}
Daley, R.: {Linear non-divergent mass-wind laws on the sphere}.
\newblock Tellus \textbf{35A}, 17--27 (1983)

\bibitem{franken2024singlelayer}
Franken, A., Caliaro, M., Cifani, P., Geurts, B.: Zeitlin truncation of a shallow water quasi-geostrophic model for planetary flow.
\newblock J. Adv. Model. Earth Syst. \textbf{16}(6), e2023MS003901 (2024)

\bibitem{franken2025casimir}
Franken, A.D., Luesink, E., Ephrati, S.R., Geurts, B.J.: Casimir preserving numerical method for global multilayer geostrophic turbulence.
\newblock J. Comput. Phys. p. 114155 (2025)

\bibitem{HaiLubWan}
Hairer, E., Lubich, C., Wanner, G.: Geometric Numerical Integration. Structure-Preserving Algorithms for Ordinary Differential Equations.
\newblock Springer-Verlag Berlin Heidelberg (2006)

\bibitem{HazMorr}
Hazeltine, R., Morrison, P.: {Hamiltonian formulation of reduced magnetohydrodynamics}.
\newblock Phys. Fluids \textbf{27}(4), 886--897 (1984)

\bibitem{HolmKuper1983}
Holm, D., Kupershmidt, B.: {Poisson brackets and {C}lebsch representations for magnetohydrodynamics, multifluid plasmas, and elasticity}.
\newblock Physica D: Nonlinear Phenomena \textbf{6}(3), 347--363 (1983)

\bibitem{HolmLuesink2021}
Holm, D., Luesink, E.: {Stochastic Wave–Current Interaction in Thermal Shallow Water Dynamics}.
\newblock J. Nonlinear Sci. \textbf{31}, 29 (2021)

\bibitem{HolmLuesinkPan2021}
Holm, D., Luesink, E., Pan, W.: {Stochastic mesoscale circulation dynamics in the thermal ocean}.
\newblock Phys. Fluids \textbf{33}, 046603 (2021)

\bibitem{HMR}
Holm, D., Marsden, J., Ratiu, T.: The {E}uler--{P}oincar\'{e} equations and semidirect products with applications to continuum theories.
\newblock Adv. Math. \textbf{137}, 1--81 (1998)

\bibitem{HoppYau}
Hoppe, J., Yau, S.T.: {Some Properties of Matrix Harmonics on $\mathrm{S}^{2}$}.
\newblock Commun. Math. Phys. \textbf{195}, 67--77 (1998)

\bibitem{KhMisMod}
Khesin, B., Misiolek, G., Modin, K.: {Geometric Hydrodynamics and Infinite-Dimensional Newton's Equations}.
\newblock Bull. Amer. Math. Soc. \textbf{58}, 377--442 (2021)

\bibitem{kraichnan1967inertial}
Kraichnan, R.H.: Inertial ranges in two-dimensional turbulence.
\newblock Phys. Fluids \textbf{10}(7), 1417 (1967)

\bibitem{Lorenz1960}
Lorenz, E.: {Energy and numerical weather prediction}.
\newblock Tellus \textbf{12}, 364--373 (1960)

\bibitem{LFEG2024}
Luesink, E., Franken, A., Ephrati, S., Geurts, B.: {Geometric derivation and structure-preserving simulation of quasi-geostrophy on the sphere}.
\newblock arXiv:2402.13707 pp. 1--16 (2024)

\bibitem{ModRoop}
Modin, K., Roop, M.: {Spatio--temporal Lie--Poisson discretization for incompressible magnetohydrodynamics on the sphere}.
\newblock IMA J. Numer. Anal. pp. 1--36 (2025)

\bibitem{ModViv1}
Modin, K., Viviani, M.: {A Casimir preserving scheme for long-time simulation of spherical ideal hydrodynamics}.
\newblock J. Fluid Mech. \textbf{884}, A22 (2020)

\bibitem{ModViv}
Modin, K., Viviani, M.: {Lie–Poisson Methods for Isospectral Flows}.
\newblock Found. Comput. Math. \textbf{20}, 889--921 (2020)

\bibitem{ModViv2}
Modin, K., Viviani, M.: {Canonical scale separation in two-dimensional incompressible hydrodynamics}.
\newblock J. Fluid Mech. \textbf{943}, A36 (2022)

\bibitem{MoGr1980}
Morrison, P., Greene, J.: {Noncanonical Hamiltonian Density Formulation of Hydrodynamics and Ideal Magnetohydrodynamics}.
\newblock Phys. Rev. Lett. \textbf{45}(10), 790--794 (1980)

\bibitem{BrienReid1967}
O'Brien, J., Reid, R.: {The Non-Linear Response of a Two-Layer, Baroclinic Ocean to a Stationary, Axially-Symmetric Hurricane: Part I. Upwelling Induced by Momentum Transfer}.
\newblock J. Atmos. Sci. \textbf{24}, 197--207 (1967)

\bibitem{Ripa1993}
Ripa, P.: {Conservation laws for primitive equations models with inhomogeneous layers}.
\newblock Geophys. Astrophys. Fluid Dyn. \textbf{70}, 85--111 (1993)

\bibitem{ripa1995low}
Ripa, P.: Low frecuency approximation of a vertically averaged ocean model with thermodynamics.
\newblock Rev. Mex. Fís. \textbf{42}(1), 117--135 (1995)

\bibitem{Ripa1995}
Ripa, P.: {On improving a one-layer ocean model with thermodynamics}.
\newblock J. Fluid Mech. \textbf{303}, 169--201 (1995)

\bibitem{Ripa1999}
Ripa, P.: {On the validity of layered models of ocean dynamics and thermodynamics with reduced vertical resolution}.
\newblock Dynam. Atmos. Ocean \textbf{29}, 1--40 (1999)

\bibitem{SchubTaftSilv2008}
Schubert, W., Taft, R., Silvers, L.: Shallow water quasi-geostrophic theory on the sphere.
\newblock J. Adv. Model. Earth Syst. \textbf{1}(2), 1--17 (2009)

\bibitem{Strauss1976}
Strauss, H.: {Nonlinear, three-dimensional magnetohydrodynamics of noncircular tokamaks}.
\newblock Phys. Fluids \textbf{19}, 134--140 (1976)

\bibitem{Verkley2009}
Verkley, W.: {A Balanced Approximation of the One-Layer Shallow-Water Equations on a Sphere}.
\newblock J. Atmos. Sci. \textbf{66}, 1735--1748 (2009)

\bibitem{WarnefordDellar2013}
Warneford, E., Dellar, P.: {The quasi-geostrophic theory of the thermal shallow water equations}.
\newblock J. Fluid. Mech. \textbf{723}, 374--403 (2013)

\bibitem{YanoMuletBonazzola2009}
Yano, J.I., Mulet, S., Bonazzola, M.: {Tropical large-scale circulations: asymptotically non-divergent?}
\newblock Tellus \textbf{61A}, 417--427 (2009)

\bibitem{Ze1991}
Zeitlin, V.: {Finite-mode analogs of {$2$}{D} ideal hydrodynamics: coadjoint orbits and local canonical structure}.
\newblock Phys. D \textbf{49}(3), 353--362 (1991)

\bibitem{Zeit}
Zeitlin, V.: {Self-consistent-mode approximation for the hydrodynamics of an incompressible fluid on non rotating and rotating spheres}.
\newblock Phys. Rev. Lett. \textbf{93}(26), 264501 (2004)

\bibitem{Ze2005}
Zeitlin, V.: {On self-consistent finite-mode approximations in (quasi-) two-dimensional hydrodynamics and magnetohydrodynamics}.
\newblock Phys. Lett. A \textbf{339}(3-5), 316--324 (2005)

\bibitem{Zeitlin2018}
Zeitlin, V.: Understanding (Almost) Everything with Rotating Shallow Water Models.
\newblock Oxford University Press (2018)

\end{thebibliography}
\end{document}